# DISTRIBUTION FREE GOODNESS-OF-FIT TESTS FOR LINEAR PROCESSES[1]

By Miguel A. Delgado, Javier Hidalgo and Carlos Velasco

*Universidad Carlos III, London School of Economics and Universidad Carlos III*

This article proposes a class of goodness-of-fit tests for the autocorrelation function of a time series process, including those exhibiting long-range dependence. Test statistics for composite hypotheses are functionals of a (approximated) martingale transformation of the Bartlett $T_p$-process with estimated parameters, which converges in distribution to the standard Brownian motion under the null hypothesis. We discuss tests of different natures such as omnibus, directional and Portmanteau-type tests. A Monte Carlo study illustrates the performance of the different tests in practice.

**1. Introduction and statement of the problem.** Let $f$ be the spectral density function of a second-order stationary time series process $\{X(t)\}_{t\in\mathbb{Z}}$ with mean $\mu$ and covariance function

$$\mathrm{Cov}(X(j),X(0)) = \int_{-\pi}^{\pi} f(\lambda)\cos(\lambda j)\,d\lambda, \qquad j=0,\pm1,\pm2,\ldots.$$

We shall assume that $\{X(t)\}_{t\in\mathbb{Z}}$ admits the Wold representation

$$X(t) = \mu + \sum_{j=0}^{\infty} a(j)\varepsilon(t-j) \qquad \text{with } a(0)=1 \text{ and } \sum_{j=0}^{\infty} a^2(j) < \infty, \tag{1}$$

for some sequence $\{\varepsilon(t)\}_{t\in\mathbb{Z}}$ satisfying $\mathbb{E}(\varepsilon(t))=0$, and $\mathbb{E}(\varepsilon(0)\varepsilon(t))=\sigma^2$ if $t=0$ and $=0$ for all $t\neq 0$. Under (1), the spectral density function of

Received May 2002; revised December 2004.

[1]Supported in part by the Spanish Dirección General de Enseñanza Superior (DGES) reference numbers BEC2001-1270 and SEJ2004-04583/ECON and by the Economic and Social Research Council (ESRC) reference number R000239936.

$\{X(t)\}_{t\in\mathbb{Z}}$ can be factorized as

$$f(\lambda) = \frac{\sigma^2}{2\pi} h(\lambda), \qquad \lambda \in [0,\pi],$$

with $h(\lambda) := |\sum_{j=0}^{\infty} a(j)e^{ij\lambda}|^2$.

Let

$$\mathcal{H} = \left\{ h_\theta : \int_0^\pi \log h_\theta(\lambda)\, d\lambda = 0, \theta \in \Theta \right\}, \tag{2}$$

where $\Theta \subset \mathbb{R}^p$ is a compact parameter space. Much of the existing time series literature is concerned with parametric estimation and testing, assuming that $h$ belongs to $\mathcal{H}$, that is, $h = h_{\theta_0}$ for some $\theta_0 \in \Theta$, because the parameter $\theta_0$ and the functional form of $h_\theta$ summarize the autocorrelation structure of $\{X(t)\}_{t\in\mathbb{Z}}$. Notice that $h \in \mathcal{H}$ in (2) guarantees that $a(0) = 1$ in (1) and $\sigma^2 = \min_{\theta\in\Theta} 2\int_0^\pi f(\lambda)/h_\theta(\lambda)\, d\lambda$. For our purposes, $\sigma^2$ can be considered a nuisance parameter, as well as the mean $\mu$.

Classical parameterizations that accommodate alternative models are the ARMA, ARFIMA, fractional noise and Bloomfield [4] exponential models (see [35] for definitions). For instance, in an ARFIMA specification, $\mathcal{H}$ consists of all functions indexed by a parameter vector $\theta = (d, \eta', \delta')'$, where $\theta \in \Theta \subset (-1/2, 1/2) \times \mathbb{R}^{p_1} \times \mathbb{R}^{p_2}$, of the form

$$h_\theta(\lambda) = \frac{1}{|1 - e^{i\lambda}|^{2d}} \left| \frac{\Xi_\eta(e^{i\lambda})}{\Phi_\delta(e^{i\lambda})} \right|^2, \qquad \lambda \in [0,\pi], \tag{3}$$

such that $\Xi_\eta$ and $\Phi_\delta$ are the moving average and autoregressive polynomials of orders $p_1$ and $p_2$, respectively, with no common zeros, all lying outside the unit circle.

Before statistical inference on the true value $\theta_0$ is made, one needs to test the hypothesis $H_0 : h \in \mathcal{H}$, which can be equivalently stated as

$$H_0 : \frac{G_{\theta_0}(\lambda)}{G_{\theta_0}(\pi)} = \frac{\lambda}{\pi} \qquad \text{for all } \lambda \in [0,\pi] \text{ and some } \theta_0 \in \Theta, \tag{4}$$

where

$$G_\theta(\lambda) := 2 \int_0^\lambda \frac{f(\bar\lambda)}{h_\theta(\bar\lambda)}\, d\bar\lambda, \qquad \lambda \in [0,\pi].$$

Under $H_0$, $G_{\theta_0}$ is the spectral distribution function of the innovation process $\{\varepsilon(t)\}_{t\in\mathbb{Z}}$ and $G_{\theta_0}(\pi) = \sigma^2$.

Given a record $\{X(t)\}_{t=1}^T$ and a consistent estimator $\theta_T$ of $\theta_0$ under $H_0$, a natural estimator of $G_{\theta_0}$ is defined as $G_{\theta_T,T}(\lambda)$, where

$$G_{\theta,T} := \frac{2\pi}{\widetilde{T}} \sum_{j=1}^{[\widetilde{T}\lambda/\pi]} \frac{I_X(\lambda_j)}{h_\theta(\lambda_j)}, \qquad \lambda \in [0,\pi]. \tag{5}$$

Here $\widetilde{T} = [T/2]$, $[z]$ being the integer part of $z$, and for a generic time series process $\{V(t)\}_{t\in\mathbb{Z}}$,

$$I_V(\lambda_j) := \frac{1}{2\pi T}\left|\sum_{t=1}^{T} V(t)e^{it\lambda_j}\right|^2, \qquad j = 1, \dots, \widetilde{T},$$

denotes the periodogram of $\{V(t)\}_{t=1}^{T}$ evaluated at the Fourier frequency $\lambda_j = 2\pi j/T$ for positive integers $j$.

The formulation of $H_0$ in (4) suggests use of Bartlett's $T_p$-process as a basis for testing $H_0$. The $T_p$-process is defined as

$$\alpha_{\theta,T}(\lambda) := \widetilde{T}^{1/2}\left[\frac{G_{\theta,T}(\lambda)}{G_{\theta,T}(\pi)} - \frac{\lambda}{\pi}\right], \qquad \lambda \in [0, \pi].$$

Notice that $\alpha_{\theta,T}$ is scale invariant and that, for $j \neq 0 \bmod(T)$, $I_V(\lambda_j)$ is mean invariant, so omission of $j = 0$ in the definition of $G_{\theta,T}$ entails mean correction. That is, $\alpha_{\theta,T}$ is independent of both $\mu$ and $\sigma^2$.

Under short-range dependence and $H_0$, we have that

$$\max_{1\le j\le \widetilde{T}} \mathbb{E}\left|\frac{I_X(\lambda_j)}{h_{\theta_0}(\lambda_j)} - I_\varepsilon(\lambda_j)\right| = o(1);$$

see [7], Theorem 10.3.1, page 346. So, it is expected that $\alpha_{\theta_0,T}$ will be asymptotically equivalent to Bartlett's $U_p$-process for $\{\varepsilon(t)\}_{t\in\mathbb{Z}}$,

$$\alpha_T^0(\lambda) := \widetilde{T}^{1/2}\left[\frac{G_T^0(\lambda)}{G_T^0(\pi)} - \frac{\lambda}{\pi}\right], \qquad \lambda \in [0, \pi],$$

with

$$G_T^0(\lambda) := \frac{2\pi}{\widetilde{T}} \sum_{j=1}^{[\widetilde{T}\lambda/\pi]} I_\varepsilon(\lambda_j), \qquad \lambda \in [0, \pi].$$

In fact, under suitable regularity conditions, we shall show below that the aforementioned equivalence also holds under long-range dependence. Observe that the $U_p$-process $\alpha_T^0$ and the $T_p$-process $\alpha_{\theta_0,T}$ are identical when $\{X(t)\}_{t\in\mathbb{Z}}$ is a white noise process.

The $U_p$-process $\alpha_T^0$ is useful for testing simple hypotheses when the innovations $\{\varepsilon(t)\}_{t=1}^{T}$ can be easily computed, as is the case when $\{X(t)\}_{t\in\mathbb{Z}}$ is an AR model. However, there are many other models of interest whose innovations $\{\varepsilon(t)\}_{t=1}^{T}$ cannot be directly computed, for example, Bloomfield's exponential model, or difficult to obtain, as in models exhibiting long-range dependence, such as ARFIMA models. In those cases, it appears computationally much simpler to use $\alpha_{\theta_0,T}$ for testing simple hypotheses.

The empirical processes $\alpha_T^0$ and $\alpha_{\theta,T}$, with fixed $\theta$, are random elements in $D[0, \pi]$, the space of right continuous functions on $[0, \pi]$ with left-hand side limits, the càdlàg space. The functional space $D[0,\pi]$ is endowed with the Skorohod metric (see, e.g., [3]) and convergence in distribution in the corresponding topology will be denoted by "$\Rightarrow$".

Under suitable regularity conditions on $\{\varepsilon(t)\}_{t\in\mathbb{Z}}$, it is well known that

$$\alpha_T^0 \Rightarrow B_\pi^1, \tag{6}$$

where $B_\pi^1$ is the standardized tied down Brownian motion at $\pi$. In terms of the standard Brownian motion $B$ on $[0,1]$, $B_\pi^1$ can be represented as

$$B_\pi^1(\lambda) = B\left(\frac{\lambda}{\pi}\right) - \frac{\lambda}{\pi}B(1), \qquad \lambda \in [0,\pi].$$

Grenander and Rosenblatt [18] proved (6) assuming that $\{\varepsilon(t)\}_{t\in\mathbb{Z}}$ is a sequence of independent and identically distributed (i.i.d.) random variables with eight bounded moments. The i.i.d. condition was relaxed by Dahlhaus [10], who assumed that $\{\varepsilon(t)\}_{t\in\mathbb{Z}}$ behaves as a martingale difference, but still assuming eight bounded moments. Recently Klüppelberg and Mikosch [27] proved (6) under i.i.d. $\{\varepsilon(t)\}_{t\in\mathbb{Z}}$, but assuming only four bounded moments. The i.i.d. requirement is relaxed by the following assumption:

A1. The innovation process $\{\varepsilon(t)\}_{t\in\mathbb{Z}}$ satisfies $\mathbb{E}(\varepsilon(t)^r|\mathcal{F}_{t-1}) = \mu_r$ with $\mu_r$ constant ($\mu_1 = 0$ and $\mu_2 = \sigma^2$) for $r = 1,\ldots,4$ and all $t = 0,\pm 1,\ldots,$ where $\mathcal{F}_t$ is the sigma algebra generated by $\{\varepsilon(s), s \le t\}$.

Assumption A1 appears to be a minimal requirement to establish a functional central limit theorem for $\alpha_T^0$, due to the quadratic nature of the periodogram.

To establish the asymptotic equivalence between $\alpha_{\theta_0,T}$ and $\alpha_T^0$, we introduce the following smoothness assumptions on $h$:

A2. (a) $h$ is a positive and continuously differentiable function on $(0,\pi]$;
(b) $|\partial \log h(\lambda)/\partial\lambda| = O(\lambda^{-1})$ as $\lambda \to 0+$.

This condition is very general and allows for a possible singularity of $h$ at $\lambda = 0$. It holds for models exhibiting long-range dependence, like ARFIMA$(p_2, d, p_1)$ models with $d \neq 0$, as can easily be checked using (3) and that $|1 - e^{i\lambda}| = |2\sin(\lambda/2)|$.

THEOREM 1. *Assuming* A1 *and* A2, *under* $H_0$, (6) *holds and*

$$\sup_{\lambda\in[0,\pi]} |\alpha_{\theta_0,T}(\lambda) - \alpha_T^0(\lambda)| = o_p(1).$$

We can relax the location of the possible singularity in $h$ at any other frequency $\lambda \neq 0$, as in [23] or, more recently, [14], or even allow for more than one singularity. However, for notational simplicity we have taken the singularity, if any, at $\lambda = 0$. If the location of the singularity is at $\lambda^0 \neq 0$, then A2 is modified to the following:

A2′. (a) $h$ is a positive and continuously differentiable function on $[0, \lambda^0) \cup (\lambda^0, \pi]$;
(b) $|\partial \log h(\lambda)/\partial\lambda| = O(|\lambda - \lambda^0|^{-1})$ as $\lambda \to \lambda^0$.

We now comment on the results of Theorem 1. This theorem indicates that $\alpha_{\theta_0,T}$ is asymptotically pivotal. One consequence is that critical regions of tests based on a continuous functional $\varphi : D[0,\pi] \mapsto \mathbb{R}$ can be easily obtained. Different functionals $\varphi$ lead to tests with different power properties. Among them are omnibus, directional and/or Portmanteau-type tests. For example, classical functionals which lead to omnibus tests are the Kolmogorov–Smirnov ($\varphi(g) = \sup_{\lambda\in[0,\pi]} |g(\lambda)|$) and the Cramér–von Mises ($\varphi(g) = \pi^{-1} \int_0^\pi g(\lambda)^2 \, d\lambda$), whereas Portmanteau tests, defined as weighted sums of squared estimated autocorrelations of the innovations, and directional tests are obtained by choosing an appropriate functional $\varphi$; see Section 3 for details.

On the other hand, in practical situations the parameters $\theta_0$ are not known and, thus, they have to be replaced by some estimate $\theta_T$. In this situation, as Theorem 2 below shows, the $T_p$-process is no longer asymptotically pivotal and, hence, the aforementioned tests are not useful for practical purposes. The unknown critical values of functionals of the $T_p$-process with estimated parameters can be approximated with the assistance of bootstrap methods. This approach has been proposed by Chen and Romano [9] and Hainz and Dahlhaus [19] for short-range models using the $U_p$-process and by Delgado and Hidalgo [11], who allow also long-range dependence models using the $T_p$-process. Alternatively, asymptotically distribution free tests can be obtained by introducing a tuning parameter that must behave in some required way as the sample size increases. Among them, the most popular one is the Portmanteau test, although it has only been justified for testing short-range models. Box and Pierce [5] showed that the partial sum of the squared residual autocorrelations of a stationary ARMA process is approximately chi-squared distributed assuming that the number of autocorrelations considered diverges to infinity with the sample size at an appropriate rate. A different approach, in the spirit of Durbin, Knott and Taylor [12] for the classical empirical process, is that in Anderson [2], who proposed to approximate the critical values of the Cramér–von Mises tests for a stationary AR model. The method considers a truncated version of the spectral representation of $\alpha_{\theta_T,T}$ with estimated orthogonal components. The number of estimated orthogonal components must suitably increase with the sample size. A similar idea was proposed by Velilla [46] for ARMA models. Finally, another alternative uses the distance between a smooth estimator of the spectral density function and its parametric estimator under $H_0$. This approach provides asymptotically distribution free tests for short-range models assuming a suitable behavior of the smoothing parameter as the sample size diverges; see, for example, Prewitt [34] and Paparoditis [33]. However, the final outcome of all these tests depends on the arbitrary choice of the tuning/smoothing parameters, for which no relevant theory is available.

This article solves some limitations of existing asymptotically pivotal tests, only justified under short-range dependence, by considering an asymptotically pivotal transformation of $\alpha_{\theta_T,T}$ related to the cusum of recursive residuals proposed by Brown, Durbin and Evans [8]. We show that our testing procedure is valid under long-range specifications. In the next section we provide regularity conditions for the weak convergence of $\alpha_{\theta_T,T}$ and its asymptotically distribution free transformation. In Section 3 we discuss the behavior of tests of a very different nature—omnibus, directional and smooth/Portmanteau—under local alternatives converging to the null at the rate $T^{-1/2}$. Section 4 reports the results of a small Monte Carlo experiment. Some final remarks are placed in Section 5. Section 6 provides lemmata with some auxiliary results, which are employed to prove, in Section 7, the main results of the paper.

**2. Tests based on a martingale transformation of the $T_p$-process with estimated parameters.** A popular estimator of $\theta_0$ is the Whittle estimator

$$\theta_T := \arg\min_{\theta\in\Theta} G_{\theta,T}(\pi), \tag{7}$$

with $G_{\theta,T}$ defined in (5). Let us define

$$\phi_\theta(\lambda) := \frac{\partial}{\partial\theta}\log h_\theta(\lambda), \qquad S_T := \frac{1}{\widetilde{T}}\sum_{j=1}^{\widetilde{T}} \phi_{\theta_0}(\lambda_j)\phi'_{\theta_0}(\lambda_j)$$

and introduce the following assumptions:

A3. (a) $\phi_{\theta_0}$ is a continuously differentiable function on $(0,\pi]$; (b) $\|\partial\phi_{\theta_0}(\lambda)/\partial\lambda\| = O(1/\lambda)$ as $\lambda\to 0+$; and for some $0<\delta<1$ and all $\lambda\in(0,\pi]$, there exists a $K<\infty$ such that (c) $\sup_{\{\theta:\|\theta-\theta_0\|\le\delta\}}\|\phi_\theta(\lambda)\|\le K|\log\lambda|$; (d)

$$\sup_{\{\theta:\|\theta-\theta_0\|\le\delta/2\}}\frac{1}{\|\theta-\theta_0\|^2}\left|\frac{h_{\theta_0}(\lambda)}{h_\theta(\lambda)} - 1 + \phi'_{\theta_0}(\lambda)(\theta-\theta_0)\right| \le \frac{K}{\lambda^\delta}\log^2\lambda;$$

and (e) $\Sigma_{\theta_0} := \pi^{-1}\int_0^\pi \phi_{\theta_0}(\lambda)\phi'_{\theta_0}(\lambda)\,d\lambda$ is positive definite.

These assumptions are standard when analyzing the asymptotic distribution of the Whittle estimator $\theta_T$ and they are satisfied for all parametric linear processes used in practice. Standard ARMA models satisfy a stronger condition, replacing the upper bounds in A3(c) and (d) by a constant independent of $\lambda$. It can easily be shown that A3 is satisfied for ARFIMA models. Note that A3(e) and Lemma 1 in Section 6 imply that $S_T$ is positive definite for $T$ large enough.

A4. The estimator in (7) satisfies the asymptotic linearization

$$\widetilde{T}^{1/2}(\theta_T - \theta_0) = S_T^{-1} \int_0^\pi \phi_{\theta_0}(\lambda)\alpha_{\theta_0,T}(d\lambda) + o_p(1). \tag{8}$$

The expansion (8), in assumption A4, is satisfied under A1–A3 and additional standard identification conditions; see [15, 20] or [45] for a later reference.

Define

$$\alpha_\infty(\lambda) := B_\pi^1(\lambda) - \left(\frac{1}{\pi}\int_0^\lambda \phi'_{\theta_0}(\bar{\lambda})\,d\bar{\lambda}\right)\Sigma_{\theta_0}^{-1}\int_0^\pi \phi_{\theta_0}(\bar{\lambda})B_\pi^1(d\bar{\lambda}).$$

THEOREM 2. *Under $H_0$ and assuming* A1–A4, *uniformly in $\lambda \in [0,\pi]$:*

(a) $$\alpha_{\theta_T,T}(\lambda) = \alpha_T^0(\lambda) - \left(\frac{1}{\widetilde{\widetilde{T}}}\sum_{j=1}^{[\widetilde{T}\lambda/\pi]} \phi'_{\theta_0}(\lambda_j)\right) S_T^{-1}\int_0^\pi \phi_{\theta_0}(\bar{\lambda})\alpha_T^0(d\bar{\lambda}) + o_p(1);$$

(b) $\alpha_{\theta_T,T} \Rightarrow \alpha_\infty$.

Theorem 2 indicates that the asymptotic critical values of tests based on $\alpha_{\theta_T,T}$ cannot be tabulated. However, we can use a transformation of $\alpha_{\theta_T,T}$ that converges in distribution to the standard Brownian motion. To this end, it is of interest to realize that Theorem 2(a) provides an asymptotic representation of $\alpha_{\theta_T,T}$ as a scaled cumulative sum (cusum) of the least squares residuals in an artificial regression model. For that purpose, observe that by (2), and using the fact that $\phi_{\theta_0}$ is integrable [A3(c)],

$$\int_0^\pi \phi_{\theta_0}(\lambda)\,d\lambda = 0. \tag{9}$$

Now, because Lemma 1 in Section 6 with $\zeta(\lambda) = \phi_{\theta_0}(\lambda)$ and (9) imply that $\|\sum_{k=1}^{\widetilde{T}} \phi_{\theta_0}(\lambda_k)\| = O(\log T)$, the uniform asymptotic expansion in Theorem 2(a) indicates that

$$\sup_{\lambda\in[0,\pi]}\left|\alpha_{\theta_T,T}(\lambda) - \frac{2\pi}{G_T^0(\pi)}\frac{1}{\widetilde{T}^{1/2}}\sum_{j=1}^{[\widetilde{T}\lambda/\pi]} u_T(j)\right| = o_p(1),$$

where

$$u_T(j) = I_\varepsilon(\lambda_j) - \gamma'_{\theta_0}(\lambda_j)\left[\sum_{k=1}^{\widetilde{T}} \gamma_{\theta_0}(\lambda_k)\gamma'_{\theta_0}(\lambda_k)\right]^{-1} \sum_{k=1}^{\widetilde{T}} \gamma_{\theta_0}(\lambda_k) I_\varepsilon(\lambda_k),$$

$$j = 1, \ldots, \widetilde{T},$$

are the least squares residuals in an artificial regression model with dependent variable $I_\varepsilon(\lambda_j)$ and a vector of explanatory variables $\gamma_{\theta_0}(\lambda_j) := (1, \phi'_{\theta_0}(\lambda_j))'$. This fact suggests employing the cusum of recursive residuals for constructing asymptotically pivotal tests, as were proposed by Brown, Durbin and Evans [8]; see also [39].

Let us define

$$A_{\theta,T}(j) := \frac{1}{\widetilde{T}} \sum_{k=j+1}^{\widetilde{T}} \gamma_\theta(\lambda_k)\gamma'_\theta(\lambda_k)$$

and assume the following:

A5. $A_{\theta_0,T}(\overline{T})$ is nonsingular for $\overline{T} = \widetilde{T} - p - 1$.

The (scaled) cusum of forward recursive least squares residuals is defined as

$$\beta_T^0(\lambda) := \frac{2\pi}{G_T^0(\pi)} \frac{1}{\widetilde{T}^{1/2}} \sum_{j=1}^{[\overline{T}\lambda/\pi]} e_T(j), \qquad \lambda \in [0, \pi],$$

where

$$e_T(j) := I_\varepsilon(\lambda_j) - \gamma'_{\theta_0}(\lambda_j) b_T(j), \qquad j = 1, \ldots, \overline{T},$$

are the forward least squares residuals and

$$b_T(j) := A_{\theta_0,T}^{-1}(j) \frac{1}{\widetilde{T}} \sum_{k=j+1}^{\widetilde{T}} \gamma_{\theta_0}(\lambda_k) I_\varepsilon(\lambda_k).$$

It is worth observing that the motivation to employ only the first $\overline{T}$ Fourier frequencies to compute the recursive residuals is due to the singularity of $A_{\theta,T}(j)$ for all $j > \overline{T}$.

The empirical process $\beta_T^0$ can be written as a linear transformation of $\alpha_T^0$,

$$\beta_T^0(\lambda) = \mathcal{L}_{\theta_0,T}\alpha_T^0(\lambda), \qquad \lambda \in [0, \pi],$$

where, for any function $g \in D[0, \pi]$,

$$\mathcal{L}_{\theta,T} g(\lambda) = g\left(\frac{\overline{T}}{\widetilde{T}}\lambda\right) - \frac{1}{\widetilde{T}} \sum_{j=1}^{[\overline{T}\lambda/\pi]} \gamma'_\theta(\lambda_j) A_{\theta,T}^{-1}(j) \int_{\lambda_{j+1}}^{\pi} \gamma_\theta(\tilde{\lambda}) g(d\tilde{\lambda}).$$

The transformation $\mathcal{L}_{\theta_0,T}$ has the limiting version $\mathcal{L}^0$, defined as

$$\mathcal{L}^0 g(\lambda) = g(\lambda) - \frac{1}{\pi}\int_0^\lambda \gamma'_{\theta_0}(\bar\lambda) A^{-1}_{\theta_0}(\bar\lambda) \int_{\bar\lambda}^\pi \gamma_{\theta_0}(\tilde\lambda) g(d\tilde\lambda)\, d\bar\lambda,$$

where

$$A_{\theta_0}(\lambda) := \int_\lambda^\pi \gamma_{\theta_0}(\tilde\lambda)\gamma'_{\theta_0}(\tilde\lambda)\, d\tilde\lambda.$$

Notice that $\mathcal{L}^0\alpha_\infty$ is the martingale innovation of $\alpha_\infty$; see [25].

This type of martingale transformation has been used by Khmaladze [25] and Aki [1] in the standard goodness-of-fit testing problem, by Nikabadze and Stute [32] for goodness-of-fit of distribution functions under random censorship, by Stute, Thies and Zhu [42], Koul and Stute [28, 29] and Khmaladze and Koul [26] for dynamic regression models, and by Stute and Zhu [43] for generalized linear models.

Henceforth, $B_\pi(\lambda) := B(\lambda/\pi)$ for $\lambda \in [0, \pi]$.

THEOREM 3. *Under* $H_0$ *and assuming* A1–A5,

$$\beta_T^0 \Rightarrow B_\pi.$$

Because $\beta_T^0$ cannot be computed in practice, as it depends on $\theta_0$, it is suggested one use $\beta_{\theta_T,T}$, where

$$\begin{aligned}\beta_{\theta,T}(\lambda) &:= \mathcal{L}_{\theta,T}\alpha_{\theta,T}(\lambda)\\ &= \frac{2\pi}{G_{\theta,T}(\pi)}\frac{1}{\widetilde{T}^{1/2}}\sum_{j=1}^{[\overline{T}\lambda/\pi]} e_{\theta,T}(j), \qquad \lambda \in [0,\pi],\end{aligned}$$

and

$$e_{\theta,T}(j) = \frac{I_X(\lambda_j)}{h_\theta(\lambda_j)} - \gamma'_\theta(\lambda_j) b_{\theta,T}(j), \qquad j = 1, \dots, \overline{T},$$

are the forward recursive residuals in the linear projection of $I_X(\lambda_j)/h_\theta(\lambda_j)$ on $\gamma_\theta(\lambda_j)$, and where

$$b_{\theta,T}(j) = A^{-1}_{\theta,T}(j)\frac{1}{\widetilde{T}}\sum_{k=j+1}^{\widetilde{T}} \gamma_\theta(\lambda_k)\frac{I_X(\lambda_k)}{h_\theta(\lambda_k)}.$$

In order to establish the asymptotic equivalence between $\beta_T^0$ and $\beta_{\theta_T,T}$, we also need some extra smoothness assumptions on the model under the null.

A6. For some $0<\delta<1$ and all $\lambda\in(0,\pi]$, there exists a constant $K<\infty$ such that

$$\sup_{\{\theta\,:\,\|\theta-\theta_0\|\le\delta\}}\frac{1}{\|\theta-\theta_0\|^2}\|\phi_\theta(\lambda)-\phi_{\theta_0}(\lambda)-\dot{\phi}_{\theta_0}(\lambda)(\theta-\theta_0)\|\le K|\log\lambda|,$$

and $\dot{\phi}_\theta$ satisfies A3(a)–(c).

This assumption holds for all models used in practice, such as ARFIMA in (3), Bloomfield's exponential model and the fractional noise models mentioned before. In fact, they satisfy even the stronger condition with $K|\log\lambda|$ replaced by $K$.

Theorem 4. *Under* $H_0$ *and assuming* A1–A6,

$$\sup_{\lambda\in[0,\pi]}|\beta_{\theta_T,T}(\lambda)-\beta_T^0(\lambda)|=o_p(1).$$

Theorem 4 holds true, mutatis mutandis, with $\theta_T$ replaced by any $T^{1/2}$-consistent estimator. Also, from a computational point of view, it is worth observing that

$$A_{\theta,T}^{-1}(j)=A_{\theta,T}^{-1}(j+1)-\frac{A_{\theta,T}^{-1}(j+1)\gamma_\theta(\lambda_j)\gamma_\theta'(\lambda_j)A_{\theta,T}^{-1}(j+1)}{\widetilde{T}+\gamma_\theta'(\lambda_j)A_{\theta,T}^{-1}(j+1)\gamma_\theta(\lambda_j)}$$

and

$$b_{\theta,T}(j)=b_{\theta,T}(j+1)+A_{\theta,T}^{-1}(j)\gamma_\theta(\lambda_j)\left[\frac{I_X(\lambda_j)}{h_\theta(\lambda_j)}-\gamma_\theta'(\lambda_j)b_{\theta,T}(j+1)\right];$$

see [8] for similar arguments.

Alternatively to $\beta_{\theta_T,T}$, we could have considered the cusum of backward recursive residuals, that is,

$$\bar{\beta}_{\theta_T,T}(\lambda):=\frac{2\pi}{G_{\theta_T,T}(\pi)}\frac{1}{\widetilde{T}^{1/2}}\sum_{j=p+1}^{[\widetilde{T}\lambda/\pi]}\bar{e}_{\theta_T,T}(j),\qquad\lambda\in[0,\pi],$$

where

$$\bar{e}_{\theta,T}(j):=\frac{I_X(\lambda_j)}{h_\theta(\lambda_j)}-\gamma_\theta'(\lambda_j)\bar{b}_{\theta,T}(j),\qquad j=p+1,\ldots,\widetilde{T},$$

$$\bar{b}_{\theta,T}(j):=\bar{A}_{\theta,T}^{-1}(j)\frac{1}{\widetilde{T}}\sum_{k=1}^{j-1}\gamma_\theta(\lambda_k)\frac{I_X(\lambda_k)}{h_\theta(\lambda_k)}\quad\text{and}\quad\bar{A}_{\theta,T}(j):=\frac{1}{\widetilde{T}}\sum_{k=1}^{j-1}\gamma_\theta(\lambda_k)\gamma_\theta'(\lambda_k).$$

In this case, we can take advantage of the computational formulae

$$\bar{A}_{\theta,T}^{-1}(j+1)=\bar{A}_{\theta,T}^{-1}(j)-\frac{\bar{A}_{\theta,T}^{-1}(j)\gamma_\theta(\lambda_{j+1})\gamma_\theta'(\lambda_{j+1})\bar{A}_{\theta,T}^{-1}(j)}{\widetilde{T}+\gamma_\theta'(\lambda_{j+1})\bar{A}_{\theta,T}^{-1}(j)\gamma_\theta(\lambda_{j+1})}$$

and

$$\bar{b}_{\theta,T}(j+1) = \bar{b}_{\theta,T}(j) + \bar{A}^{-1}_{\theta,T}(j+1)\gamma_\theta(\lambda_{j+1})\left[\frac{I_X(\lambda_{j+1})}{h_\theta(\lambda_{j+1})} - \gamma'_\theta(\lambda_{j+1})\bar{b}_{\theta,T}(j)\right].$$

This formulation may be useful in small samples when we suspect that the main discrepancy between the null and the alternative is near $\pi$. However, from Theorems 3 and 4, it is easily seen that the empirical processes $\bar{\beta}_{\theta_T,T}$ and $\beta_{\theta_T,T}$ have the same asymptotic behavior.

Let $\varphi : D[0,\pi] \to \mathbb{R}$ be a continuous functional. Under $H_0$ and the conditions in Theorem 4,

$$\varphi(\beta_{\theta_T,T}) \xrightarrow{d} \varphi(B_\pi),$$

as a consequence of the continuous mapping theorem. For instance,

$$\hat{K}_T = \sup_{j=1,\dots,\overline{T}} \left|\beta_{\theta_T,T}\left(\frac{j\pi}{\overline{T}}\right)\right| \xrightarrow{d} \sup_{\lambda\in[0,\pi]} |B_\pi(\lambda)| \stackrel{d}{=} \sup_{\omega\in[0,1]} |B(\omega)|,$$

$$\hat{C}_T = \frac{1}{\overline{T}}\sum_{j=1}^{\overline{T}} \beta_{\theta_T,T}\left(\frac{j\pi}{\overline{T}}\right)^2 \xrightarrow{d} \frac{1}{\pi}\int_0^\pi B_\pi^2(\lambda)\,d\lambda \stackrel{d}{=} \int_0^1 B^2(\omega)\,d\omega.$$

The above limiting distributions are tabulated; see, for example, [40], pages 34 and 748.

**3. Local alternatives: omnibus, directional and Portmanteau tests.** In this section we shall show that tests based on $\beta_{\theta_T,T}$ are able to detect local alternatives of the type

$$H_{1T} : h(\lambda) = h_{\theta_0}(\lambda)\left(1 + \tau\frac{1}{\widetilde{T}^{1/2}}l(\lambda) + \frac{1}{\widetilde{T}}s_T(\lambda)\right),$$

$$\lambda \in [0,\pi] \text{ and for some } \theta_0 \in \Theta,$$

where $\int_0^\pi l(\lambda)\,d\lambda = 0$, $l(\lambda)$ satisfies the same properties as $\phi_{\theta_0}$ in A3(a)–(c), $\tau$ is a constant, possibly unknown, and for some finite $T_0$, $|s_T(\cdot)|$ is integrable for all $T > T_0$. Let us consider some examples.

Example 1. If we wish to study departures from the white noise hypothesis in the direction of fractional alternatives, we have

$$\frac{h(\lambda)}{h_{\theta_0}(\lambda)} = \frac{1}{|2\sin(\lambda/2)|^{2d/\widetilde{T}^{1/2}}}, \qquad \lambda \in [0,\pi],$$

for some $d \neq 0$. By a simple Taylor expansion up to the second term,

$$l(\lambda) = -2\log|2\sin(\lambda/2)| \quad \text{and} \quad \tau = d,$$

respectively, with the remainder function $s_T$ being such that, for some $0 \le \epsilon < 1$, $|s_T(\lambda)| \le K|\lambda|^{-\epsilon}$ for all large $T$ and some $K < \infty$.

EXAMPLE 2. If we consider departures in the direction of MA(1) alternatives, we obtain

$$\frac{h(\lambda)}{h_{\theta_0}(\lambda)} = 1 - \eta \frac{1}{\widetilde{T}^{1/2}} 2\cos(\lambda) + \frac{1}{\widetilde{T}}\eta^2, \qquad \lambda \in [0, \pi].$$

Thus, $\tau = \eta$, $l(\lambda) = -2\cos(\lambda)$ and $s_T(\lambda) = \eta^2$.

EXAMPLE 3. If we consider departures in the direction of AR(1) alternatives, then

$$\frac{h(\lambda)}{h_{\theta_0}(\lambda)} = \left[1 - \delta \frac{1}{\widetilde{T}^{1/2}} 2\cos(\lambda) + \frac{1}{\widetilde{T}}\delta^2\right]^{-1}, \qquad \lambda \in [0, \pi].$$

Thus, $\tau = \delta$ and $l(\lambda) = 2\cos(\lambda)$ with $|s_T(\lambda)| \le K$, for all large $T$ and some $K < \infty$.

For $\lambda \in [0, \pi]$, let us define

$$L(\lambda) := \frac{1}{\pi}\int_0^\lambda \left\{ l(\bar{\lambda}) - \gamma'_{\theta_0}(\bar{\lambda}) A_{\theta_0}^{-1}(\bar{\lambda}) \frac{1}{\pi}\int_{\bar{\lambda}}^\pi \gamma_{\theta_0}(\tilde{\lambda}) l(\tilde{\lambda})\, d\tilde{\lambda} \right\} d\bar{\lambda} \tag{10}$$

and

$$M(\lambda) := B_\pi(\lambda) + \tau \cdot L(\lambda), \qquad \lambda \in [0, \pi].$$

We have the following theorem.

THEOREM 5. *Assuming the same assumptions as in Theorem 4, under $H_{1T}$,*

$$\beta_{\theta_T, T} \Rightarrow M.$$

Using the fact that $M$ and $B_\pi$ are identically distributed, except for the deterministic shift $\tau \cdot L$, and taking into account that $2^{1/2}\sin((j - 1/2)\lambda)$ and $1/(j-1/2)^2\pi^2$ are the eigenfunctions and eigenvalues in the Kac–Siegert representation of $B_\pi$ [24], the orthogonal components of $M$,

$$m(j) := 2^{1/2}(j - \tfrac{1}{2}) \int_0^\pi \sin((j - \tfrac{1}{2})\lambda) M(\lambda)\, d\lambda, \qquad j = 1, 2, \ldots,$$

are independently distributed normal random variables with mean $\tau \cdot \ell(j)$ and variance 1, where

$$\ell(j) = 2^{1/2}(j - \tfrac{1}{2}) \int_0^\pi \sin((j - \tfrac{1}{2})\lambda) L(\lambda)\, d\lambda, \qquad j = 1, 2, \ldots.$$

Using the (asymptotically) orthogonal components of $\beta_{\theta_T, T}$,

$$\tilde{m}_T(j) = 2^{1/2}(j - \tfrac{1}{2}) \int_0^\pi \sin((j - \tfrac{1}{2})\lambda) \beta_{\theta_T, T}(\lambda)\, d\lambda, \qquad j = 1, 2, \ldots,$$

we obtain the spectral representation,

$$\beta_{\theta_T,T}(\lambda) = 2^{1/2} \sum_{j=1}^{\infty} \frac{\tilde{m}_T(j) \sin((j-1/2)\lambda)}{\pi(j-1/2)}, \qquad \lambda \in [0,\pi].$$

By Theorem 5 and the continuous mapping theorem, finitely many of the $\tilde{m}_T(j)$'s converge in distribution to the corresponding $m(j)$'s under $H_{1T}$. Using Parseval's theorem,

$$\hat{C}_T \xrightarrow{d} \sum_{j=1}^{\infty} \frac{m^2(j)}{(j-1/2)^2 \pi^2}.$$

Using similar arguments to those in [13] in the context of the standard empirical process with estimated parameters, tests based on

$$\tilde{W}_{n,T} := \sum_{j=1}^{n} \tilde{m}_T^2(j),$$

with a reasonable choice of $n \geq 1$, will lead to gains in power, compared to $\hat{C}_T$, in the direction of alternatives with significant autocorrelations at high lags. These Portmanteau tests are related to Neyman's [31] smooth tests, a compromise between omnibus and directional tests, and for each $n \geq 1$, under $H_{1T}$ we have that

$$\tilde{W}_{n,T} \xrightarrow{d} \chi_n^2\left(\tau^2 \sum_{j=1}^{n} \ell^2(j)\right).$$

That is, tests based on $\tilde{W}_{n,T}$ are asymptotically pivotal under $H_0$ ($\tau = 0$) for each choice of $n$, and more importantly, they are able to detect local alternatives converging to the null at the parametric rate $T^{-1/2}$, provided that $\ell(j) \neq 0$ for some $j = 1, \ldots, n$. The latter is in contrast with the classical Portmanteau tests based on

$$\tilde{Q}_{n_T,T} := \sum_{j=1}^{n_T} (T^{1/2} \tilde{\rho}_T(j))^2, \tag{11}$$

where $\tilde{\rho}_T(j)$ is some estimate of the $j$th autocorrelation of the residuals. It has been shown that $\tilde{Q}_{n_T,T}$ is approximately distributed as a $\chi^2_{n_T - p}$ under $H_0$ specifying a short-range model and assuming that $n_T$ diverges as $T \to \infty$. On the other hand, the resulting test is able to detect alternatives converging to the null at the rate $n_T^{1/4} T^{-1/2}$ (see, e.g., [21]), which is slower than $T^{-1/2}$.

In practice, it is recommended that one use the discrete version

$$\hat{W}_{n,T} := \sum_{j=1}^{n} \hat{m}_T^2(j)$$

of $\tilde{W}_{n,T}$, with

$$\hat{m}_T(j) := 2^{1/2}\left(j-\frac{1}{2}\right)\cdot\frac{\pi}{\overline{\overline{T}}}\sum_{k=1}^{\overline{T}}\sin\left(\left(j-\frac{1}{2}\right)\frac{\pi k}{\overline{\overline{T}}}\right)\beta_{\theta_T,T}\left(\frac{\pi k}{\overline{\overline{T}}}\right).$$

On the other hand, optimal tests of $H_0$ in the direction $H_{1T}$ can be constructed applying results in [16] (see also [17] and references therein), as was suggested by Stute [41] in the context of goodness-of-fit testing of a regression function. Asymptotically, testing for $H_0$ in the direction of $H_{1T}$ is equivalent to testing $\bar{H}_0 : \mathbb{E}(m(j)) = 0$ for all $j \geq 1$, against $\bar{H}_1 : \mathbb{E}(m(j)) = \tau \cdot \ell(j)$ for all $j \geq 1$ with $L$ known, but maybe with unknown $\tau$. Under $\bar{H}_0$, the distribution of $\{m(j)\}_{j\geq 1}$ is completely specified, as it is also under $\bar{H}_1$ when the parameter $\tau$ is known. Then the likelihood-ratio for a finite-dimensional set $(m(1),\ldots,m(n))$ is

$$\Lambda_n = \exp\left(\tau\sum_{j=1}^{n}\ell(j)\cdot\left(m(j)-\frac{\tau\cdot\ell(j)}{2}\right)\right). \tag{12}$$

Grenander [16] showed that $\Lambda_n \to_p \Lambda_\infty$ as $n \to \infty$, and that the most powerful test at significance level $\alpha$ has a critical region of the form $\{\Lambda_\infty > k\}$, with $P_0\{\Lambda_\infty > k\} = \alpha$ if $\sum_{j=1}^{\infty}\ell^2(j) < \infty$. The latter condition is satisfied in our context by Parseval's theorem and A3(c) because $l$ is a square integrable function.

Define

$$\psi := \frac{\sum_{j=1}^{\infty}\ell(j)\cdot m(j)}{(\sum_{j=1}^{\infty}\ell^2(j))^{1/2}}.$$

Then under $H_0$, $\psi \overset{d}{=} N(0,1)$, and in view of (12), $\psi$ forms a basis to obtain optimal critical regions. When the sign of $\tau$ is known, the critical region of the uniformly most powerful test at significance level $\alpha$ is $\{\psi > z_{1-\alpha}\}$ when $\tau > 0$ and $\{\psi < -z_{1-\alpha}\}$ when $\tau < 0$, where $z_\upsilon$ is the $\upsilon$ quantile of the standard normal. Also, when the sign of $\tau$ is unknown, the most powerful unbiased test at significance level $\alpha$ has critical region given by $\{|\psi| > z_{1-\alpha/2}\}$.

These arguments suggest an (asymptotically) optimal Neyman–Pearson test in the direction of $H_{1T}$ based on the first $n$ orthogonal components of $\beta_{\theta_T,T}$, using the test statistic

$$\hat{\psi}_{n,T} = \frac{\sum_{j=1}^{n}\ell(j)\cdot\hat{m}_T(j)}{(\sum_{j=1}^{n}\ell^2(j))^{1/2}}.$$

Schoenfeld [38] proposes the same type of statistic in the standard goodness-of-fit testing context. Under $H_0$ and the assumptions in previous sections, we have that

$$\hat{\psi}_{n,T} \overset{d}{\to} N(0,1) \qquad \text{as } T \to \infty \text{ for each fixed } n.$$

Also, arguing as in Schoenfeld's [38] Theorem 3, the convergence in distribution of $\hat{\psi}_{n_T,T}$ when $n_T$ increases with $T$ can be shown. Approximately optimal tests for $H_0$ in the direction of $H_{1T}$ reject $H_0$ at significance level $\alpha$ when $|\hat{\psi}_{n_T,T}| > z_{1-\alpha/2}$ if $\tau$ has unknown sign, $\hat{\psi}_{n_T,T} > z_{1-\alpha}$ when $\tau > 0$ and $\hat{\psi}_{n_T,T} < -z_{1-\alpha}$ when $\tau < 0$.

**4. Some Monte Carlo experiments.** A small Monte Carlo study was carried out to investigate the finite sample performance of the different tests. To that end, we considered the AR(1), MA(1) and ARFIMA$(0, d_0, 0)$ models

$$(1-\delta_0 L)X(t) = \varepsilon(t), \tag{13}$$

$$X(t) = (1-\eta_0 L)\varepsilon(t), \tag{14}$$

$$(1-L)^{d_0}X(t) = \varepsilon(t), \tag{15}$$

respectively, where the parameter $\theta_0$ equals $\delta_0$, $\eta_0$ and $d_0$ for the different models and $L$ is the lag operator. The innovations $\{\varepsilon(t)\}_{t=1}^{T}$ are i.i.d. $\mathcal{N}(0,1)$, and the sample sizes used are $T = 200$ and $500$ with different values of the parameters $\delta_0$, $\eta_0$ and $d_0$. For models (13) and (14), we considered $\delta_0$, $\eta_0 = -0.8, -0.5, 0.0, 0.5, 0.8$, whereas for model (15) $d_0 = 0.0, 0.2, 0.4$. The ARFIMA model was simulated using an algorithm by Hosking [22].

For the three models and all values of $\theta_0$, we computed the proportion of rejections in 50,000 generated samples for both sample sizes. Whittle estimates are obtained according to (7). For each of the models considered $\phi_\theta$ is given by

$$\text{AR(1)}, \qquad \theta = \delta: \phi_\delta(\lambda) = \frac{\partial}{\partial\delta}\log|1-\delta e^{i\lambda}|^{-2} = -2\frac{\delta - \cos\lambda}{1-2\delta\cos\lambda+\delta^2};$$

$$\text{MA(1)}, \qquad \theta = \eta: \phi_\eta(\lambda) = \frac{\partial}{\partial\eta}\log|1-\eta e^{i\lambda}|^{2} = 2\frac{\eta - \cos\lambda}{1-2\eta\cos\lambda+\eta^2};$$

$$\text{ARFIMA}(0,d,0), \qquad \theta = d: \phi_d(\lambda) = \frac{\partial}{\partial d}\log|1-e^{i\lambda}|^{-2d} = -2\log|2\sin(\lambda/2)|.$$

We also report, as a benchmark, the proportion of rejections using

$$C_T^0 := \frac{1}{\pi}\int_0^\pi \alpha_{\theta_0,T}^2(\lambda)\,d\lambda = T\sum_{j=1}^{\infty}\frac{\rho_{\theta_0,T}^2(j)}{\pi^2 j^2},$$

which is suitable for testing simple hypotheses. In addition, for the sake of comparison, we provide the results for the Box and Pierce [5] test statistic (11) using several values of $n_T$ increasing with $T$, where $\tilde{\rho}_T(j)$, $j \geq 1$, are the sample autocorrelations of the residuals $\{\hat{\varepsilon}(t)\}_{t=1}^{T}$. Specifically, for the AR(1) model,

$$\hat{\varepsilon}(t) = (1-\delta_T L)X(t),$$

with $X(t) = 0$ for $t \leq 0$; for the MA(1) model,

$$\hat{\varepsilon}(t) = X(t) - \eta_T \hat{\varepsilon}(t-1),$$

with $\hat{\varepsilon}(0) = 0$; and for the ARFIMA$(0, d, 0)$ model,

$$\hat{\varepsilon}(t) = \sum_{j=0}^{t-1} \vartheta(j, d_T) X(t-j),$$

where $\vartheta(j, d)$ are the coefficients in the formal expansion

$$(1-L)^d = \sum_{j=0}^{\infty} \vartheta(j, d) L^j,$$

with

$$\vartheta(j, d) = \frac{\Gamma(j-d)}{\Gamma(-d)\Gamma(j+1)}, \qquad \Gamma(a) = \int_0^{\infty} x^{a-1} e^{-x}\, dx.$$

The standardized values of $\tilde{Q}_{n_T,T}$, $(\tilde{Q}_{n_T,T} - n_T)/\sqrt{2n_T}$ are compared with the 5% critical value of the standard normal (see Hong [21]) instead of the usual $\chi^2_{(n_T-1)}$ approximation correcting by the loss of degrees of freedom due to parameter estimation, which is justified under Gaussianity. The two approximations provide a similar proportion of rejections. We also tried the weighting suggested by Ljung and Box [30], which produced very similar results.

First we analyze the size accuracy of the Cramér–von Mises test based on $\beta_{\theta_T,T}$. The empirical sizes of the tests based on $\hat{C}_T$, reported in Table 1, are reasonably close to the nominal ones. The asymptotic approximation improves noticeably when the sample size increases from $T = 200$ to $T = 500$, this improvement being uniform for all the models, although the empirical size is smaller than the nominal level. Tests based on $\tilde{Q}_{n_T,T}$ have serious size distortions for the smaller sample size and large values of $|\eta|$ in the MA(1) model, since Whittle estimates can be quite biased in these cases. The empirical size of tests based on $\tilde{Q}_{n_T,T}$ depends substantially on the number of autocorrelations used. In addition, for the larger choices of $n_T$ implemented, $\tilde{Q}_{n_T,T}$ over-rejects $H_0$. The usual recommendation $n_T = o(T^{1/2})$ also seems reasonable here, in terms of size accuracy.

Next we study the power performance of the tests. To this end, we report first, in Table 2, the proportion of rejections under the alternative hypothesis for different nonnested specifications with the model specified under the null. We cannot conclude that one test is clearly superior to the others in any of the four cases analyzed. As expected, the power of the Portmanteau test decreases as $n_T$ increases. In view of Tables 1 and 2, we can conclude that a choice of large $n_T$, around $T^{-1/2}$, produces reasonable size accuracy, but such a choice is not the best possible one in order to maximize the power. The test based on $\hat{C}_T$ is fairly powerful compared to the Portmanteau test for all cases considered, and it works remarkably well when testing an AR(1) in the direction of an MA(1) alternative.

TABLE 1
*Empirical size of omnibus and Portmanteau tests at 5% significance*

| | $T=200$ | | | | | | $T=500$ | | | | | |
|---|---|---|---|---|---|---|---|---|---|---|---|---|
| | $\hat{C}_T$ | $C_T^0$ | $\tilde{Q}_{3,T}$ | $\tilde{Q}_{6,T}$ | $\tilde{Q}_{10,T}$ | $\tilde{Q}_{20,T}$ | $\hat{C}_T$ | $C_T^0$ | $\tilde{Q}_{3,T}$ | $\tilde{Q}_{6,T}$ | $\tilde{Q}_{15,T}$ | $\tilde{Q}_{35,T}$ |
| $\delta_0, H_0:\text{AR}(1)$ | | | | | | | | | | | | |
| −0.8 | 4.92 | 4.69 | 3.34 | 3.72 | 3.91 | 3.61 | 5.07 | 5.17 | 3.56 | 3.87 | 4.35 | 3.97 |
| −0.5 | 4.38 | 4.96 | 2.80 | 3.38 | 3.60 | 3.41 | 4.96 | 5.16 | 3.12 | 3.75 | 4.17 | 3.82 |
| 0.0 | 4.07 | 4.96 | 2.66 | 3.35 | 3.45 | 3.37 | 4.62 | 5.10 | 3.00 | 3.63 | 4.11 | 3.82 |
| 0.5 | 3.59 | 4.95 | 2.67 | 3.33 | 3.57 | 3.40 | 4.50 | 5.04 | 2.97 | 3.82 | 4.17 | 3.80 |
| 0.8 | 3.08 | 4.92 | 2.89 | 3.44 | 3.73 | 3.54 | 4.27 | 5.11 | 3.33 | 3.77 | 4.32 | 3.88 |
| $\eta_0, H_0:\text{MA}(1)$ | | | | | | | | | | | | |
| −0.8 | 4.25 | 8.37 | 4.32 | 4.54 | 4.42 | 3.95 | 4.89 | 6.67 | 4.13 | 4.39 | 4.56 | 4.07 |
| −0.5 | 4.16 | 5.06 | 2.83 | 3.41 | 3.65 | 3.38 | 4.89 | 5.18 | 3.13 | 3.76 | 4.15 | 3.83 |
| 0.0 | 4.08 | 4.96 | 2.51 | 3.26 | 3.46 | 3.32 | 4.62 | 5.10 | 2.94 | 3.61 | 4.05 | 3.82 |
| 0.5 | 3.60 | 5.08 | 2.65 | 3.30 | 3.55 | 3.41 | 4.49 | 5.15 | 2.96 | 3.77 | 4.13 | 3.82 |
| 0.8 | 3.89 | 7.72 | 15.33 | 15.30 | 15.33 | 15.05 | 4.63 | 6.42 | 8.03 | 8.44 | 8.68 | 8.17 |
| $d_0, H_0:\text{I}(d)$ | | | | | | | | | | | | |
| 0.0 | 3.53 | 4.96 | 2.76 | 3.40 | 3.68 | 3.47 | 4.48 | 5.10 | 3.13 | 3.90 | 4.29 | 3.83 |
| 0.2 | 3.54 | 4.95 | 2.76 | 3.39 | 3.63 | 3.46 | 4.54 | 5.15 | 3.14 | 3.89 | 4.27 | 3.81 |
| 0.4 | 3.58 | 5.21 | 2.79 | 3.39 | 3.59 | 3.44 | 4.58 | 5.37 | 3.14 | 3.88 | 4.27 | 3.80 |

Finally, we analyze the power of the different tests when testing an AR(1) specification in the direction of local ARFIMA$(1,d,0)$ with $d=\tau/T^{1/2}$, and in the direction of local ARMA$(1,1)$ alternatives with moving average parameter $\eta=\tau/T^{1/2}$, for different values of $\tau$. The proportion of rejections for these designs is reported in Tables 3 and 4. We also consider tests based on the test statistics $\hat{W}_{n,T}$ and $\hat{\psi}_{n,T}$ (one-sided and two-sided, $\hat{\psi}^+_{n,T}$ and $|\hat{\psi}_{n,T}|$ resp.), choosing $n=3$ and 6, which has been recommended by Stute, Thies and Zhu [42] for a different goodness-of-fit test problem. Of course, tests based on the first $n$ (asymptotic) orthogonal components of $\beta_{\theta_T,T}$ are sensitive to the choice of $n$, as also happens with tests based on the $n$ (asymptotic) orthogonal components of $\alpha_{\theta_T,T}$ (the estimated autocorrelations of the innovations) in Portmanteau tests. The omnibus test based on $\hat{C}_T$ still works fairly well compared to the others, including the optimal and smooth tests. The directional tests are the most powerful in the directions for which they are designed, and the tests based on $\hat{W}_{n,T}$ and $\tilde{Q}_{n_T,T}$ work very similarly, though $\hat{W}_{n,T}$ exhibits a better size precision for the choices of $n$ considered.

Table 2

*Empirical power of omnibus and Portmanteau tests at 5% significance*

| | $T=200$ | | | | | $T=500$ | | | | |
|---|---|---|---|---|---|---|---|---|---|---|
| | $\hat{C}_T$ | $\tilde{Q}_{3,T}$ | $\tilde{Q}_{6,T}$ | $\tilde{Q}_{10,T}$ | $\tilde{Q}_{20,T}$ | $\hat{C}_T$ | $\tilde{Q}_{3,T}$ | $\tilde{Q}_{6,T}$ | $\tilde{Q}_{15,T}$ | $\tilde{Q}_{35,T}$ |
| $\eta$, $H_0$: AR(1), $H_1$: MA(1) | | | | | | | | | | |
| $-0.8$ | 100.00 | 99.97 | 99.95 | 99.25 | 92.34 | 100.00 | 100.00 | 100.00 | 100.00 | 100.00 |
| $-0.5$ | 80.82 | 70.16 | 55.53 | 44.38 | 31.25 | 99.84 | 99.23 | 97.54 | 88.65 | 68.72 |
| 0.2 | 7.12 | 5.04 | 4.98 | 4.86 | 4.34 | 12.16 | 8.31 | 7.35 | 6.27 | 5.21 |
| 0.5 | 70.82 | 72.03 | 57.50 | 46.06 | 32.15 | 98.59 | 99.32 | 97.83 | 89.19 | 69.29 |
| 0.8 | 99.56 | 99.99 | 99.95 | 99.30 | 92.76 | 100.00 | 100.00 | 100.00 | 100.00 | 100.00 |
| $\delta$, $H_0$: MA(1), $H_1$: AR(1) | | | | | | | | | | |
| $-0.8$ | 100.00 | 100.00 | 100.00 | 100.00 | 99.99 | 100.00 | 100.00 | 100.00 | 100.00 | 100.00 |
| $-0.5$ | 84.36 | 77.15 | 66.51 | 57.37 | 44.02 | 99.73 | 99.47 | 98.45 | 94.26 | 82.89 |
| 0.2 | 7.16 | 3.71 | 3.99 | 3.94 | 3.63 | 12.04 | 6.65 | 6.42 | 5.73 | 4.80 |
| 0.5 | 77.08 | 74.86 | 64.04 | 54.79 | 31.78 | 99.19 | 99.41 | 98.35 | 93.77 | 82.04 |
| 0.8 | 100.00 | 100.00 | 100.00 | 100.00 | 99.97 | 100.00 | 100.00 | 100.00 | 100.00 | 100.00 |
| $\delta$, $H_0$: I($d$), $H_1$: AR(1) | | | | | | | | | | |
| 0.2 | 11.34 | 12.84 | 13.00 | 11.27 | 13.13 | 34.92 | 33.35 | 33.01 | 23.98 | 15.71 |
| 0.5 | 26.81 | 34.11 | 41.17 | 35.55 | 24.94 | 75.29 | 81.36 | 87.81 | 80.73 | 58.52 |
| 0.8 | 9.82 | 12.86 | 21.01 | 21.32 | 15.41 | 33.21 | 38.74 | 57.53 | 61.63 | 39.15 |
| $d$, $H_0$: AR(1), $H_1$: I($d$) | | | | | | | | | | |
| 0.1 | 8.22 | 4.98 | 5.66 | 5.11 | 4.83 | 16.79 | 12.07 | 14.09 | 12.34 | 9.10 |
| 0.2 | 19.90 | 13.74 | 16.20 | 15.23 | 11.81 | 51.77 | 45.04 | 53.29 | 47.54 | 36.11 |
| 0.3 | 36.03 | 25.92 | 32.00 | 30.50 | 24.35 | 82.80 | 74.84 | 85.12 | 81.44 | 69.62 |
| 0.4 | 48.83 | 34.86 | 43.78 | 43.31 | 35.48 | 94.40 | 87.30 | 95.56 | 94.31 | 87.38 |

**5. Final remarks.** Our results can be extended to goodness-of-fit tests of models that can accommodate simultaneously stationary and nonstationary time series. For instance, if the increments $Y(t) := (1-L)X(t)$, $t = 0, \pm 1, \ldots,$ are second order stationary with zero mean and spectral density $g$ such that

$$\lim_{\lambda \to 0+} |\lambda|^{2(d-1)} g(\lambda) = G > 0 \qquad \text{for some } d \in [0.5, 1.5),$$

we can define the pseudo-spectral density function of $\{X(t)\}_{t\in\mathbb{Z}}$, $f$, as

$$f(\lambda) = \frac{1}{|1-e^{i\lambda}|^2} g(\lambda).$$

Thus, when $d \neq 1$, $g$ has a singularity at $\lambda = 0$, as happens with many long-range dependent time series (cf. A2). If $\{X(t)\}_{t\in\mathbb{Z}}$ is stationary, $f$ becomes the standard spectral density function.

If either $\{Y(t)\}_{t\in\mathbb{Z}}$ or $\{X(t)\}_{t\in\mathbb{Z}}$ satisfies Wold's decomposition, $f$ admits the factorization

$$f(\lambda) = \frac{\sigma^2}{2\pi} h(\lambda),$$

TABLE 3

*Empirical size and power under local alternatives at 5% significance*

| $\tau$ | $\rho$ | $\hat{C}_T$ | $\hat{W}_{3,T}$ | $\hat{W}_{6,T}$ | $\lvert\hat{\psi}_{3,T}\rvert$ | $\lvert\hat{\psi}_{6,T}\rvert$ | $\hat{\psi}^+_{3,T}$ | $\hat{\psi}^+_{6,T}$ | $\tilde{Q}_{3,T}$ | $\tilde{Q}_{6,T}$ |
|---|---|---|---|---|---|---|---|---|---|---|
| $H_0: \mathrm{AR}(1),\ H_1: \mathrm{ARFIMA}(1, d = \tau/T^{1/2}, 0)$ | | | | | | | | | | |
| | | | | | | $T = 200$ | | | | |
| 0 | 0.0 | 4.07 | 3.19 | 2.59 | 4.70 | 4.81 | 4.48 | 5.12 | 2.66 | 3.35 |
| | 0.5 | 3.59 | 2.98 | 2.32 | 3.79 | 4.24 | 3.62 | 3.99 | 2.67 | 3.33 |
| | 0.8 | 3.08 | 2.52 | 1.94 | 3.94 | 3.10 | 3.75 | 4.02 | 2.89 | 3.44 |
| 1 | 0.0 | 6.26 | 5.40 | 4.37 | 8.39 | 11.13 | 13.44 | 16.63 | 3.68 | 4.25 |
| | 0.5 | 3.57 | 2.90 | 2.26 | 3.45 | 4.19 | 4.19 | 5.64 | 2.73 | 3.37 |
| | 0.8 | 3.01 | 2.25 | 1.66 | 4.10 | 4.52 | 7.80 | 8.53 | 3.87 | 4.41 |
| 2 | 0.0 | 12.19 | 12.04 | 10.53 | 19.93 | 26.15 | 28.94 | 35.10 | 7.80 | 9.13 |
| | 0.5 | 3.44 | 2.91 | 2.36 | 3.47 | 4.15 | 4.25 | 6.27 | 2.91 | 3.58 |
| | 0.8 | 4.84 | 3.16 | 2.19 | 9.17 | 10.33 | 16.59 | 17.98 | 8.45 | 7.58 |
| 3 | 0.0 | 21.92 | 23.63 | 21.27 | 35.77 | 44.37 | 47.20 | 54.61 | 15.17 | 18.02 |
| | 0.5 | 3.26 | 2.74 | 2.39 | 3.65 | 4.43 | 4.99 | 6.48 | 3.27 | 3.92 |
| | 0.8 | 9.13 | 6.61 | 4.10 | 20.13 | 22.90 | 31.95 | 35.14 | 21.18 | 16.12 |
| 4 | 0.0 | 33.38 | 27.13 | 24.15 | 50.40 | 59.39 | 62.18 | 69.12 | 23.88 | 29.88 |
| | 0.5 | 3.41 | 2.47 | 2.38 | 4.09 | 4.75 | 6.80 | 7.61 | 4.32 | 4.67 |
| | 0.8 | 17.48 | 14.65 | 9.09 | 38.10 | 43.37 | 53.13 | 57.56 | 46.00 | 33.97 |
| | | | | | | $T = 500$ | | | | |
| 0 | 0.0 | 4.62 | 4.22 | 3.66 | 4.81 | 4.78 | 4.57 | 5.06 | 3.00 | 3.63 |
| | 0.5 | 4.50 | 3.99 | 3.40 | 4.26 | 4.58 | 4.27 | 4.43 | 2.97 | 3.82 |
| | 0.8 | 4.27 | 3.56 | 3.09 | 3.90 | 3.85 | 4.63 | 3.63 | 3.33 | 3.77 |
| 1 | 0.0 | 6.93 | 7.03 | 6.29 | 9.35 | 11.62 | 14.63 | 17.54 | 4.37 | 5.13 |
| | 0.5 | 4.58 | 4.42 | 4.08 | 4.85 | 5.35 | 58.30 | 7.43 | 3.02 | 3.93 |
| | 0.8 | 4.74 | 4.13 | 3.47 | 5.72 | 5.90 | 9.61 | 9.83 | 4.12 | 4.64 |
| 2 | 0.0 | 14.22 | 15.51 | 14.23 | 23.43 | 29.37 | 33.47 | 39.37 | 10.03 | 11.60 |
| | 0.5 | 4.69 | 4.72 | 4.67 | 4.83 | 6.49 | 6.37 | 10.18 | 3.08 | 4.21 |
| | 0.8 | 7.36 | 6.13 | 4.73 | 11.57 | 12.08 | 19.11 | 19.81 | 7.27 | 7.38 |
| 3 | 0.0 | 26.86 | 31.03 | 29.55 | 44.70 | 53.35 | 56.44 | 63.59 | 21.28 | 24.91 |
| | 0.5 | 4.65 | 5.04 | 5.48 | 4.71 | 7.14 | 5.44 | 11.31 | 3.30 | 4.60 |
| | 0.8 | 13.56 | 11.62 | 8.18 | 23.46 | 24.65 | 34.56 | 35.78 | 15.23 | 13.51 |
| 4 | 0.0 | 43.62 | 51.19 | 49.81 | 66.34 | 74.28 | 75.93 | 81.84 | 37.13 | 43.93 |
| | 0.5 | 4.65 | 5.18 | 6.35 | 5.05 | 7.03 | 5.09 | 10.80 | 3.81 | 5.09 |
| | 0.8 | 24.44 | 23.10 | 16.17 | 42.07 | 44.05 | 54.86 | 56.23 | 31.28 | 25.74 |

$\lvert\hat{\psi}_{n,T}\rvert$ denotes two-sided tests, whereas $\hat{\psi}^+_{n,T}$ are one-sided (right-hand side) tests.

where $h$ satisfies A2. Thus, given a parametric family $\mathcal{H}$, for example, the ARFIMA specification given in (3), a $T_p$-process for testing that $h \in \mathcal{H}$ is

$$\alpha^w_{\theta_T,T}(\lambda) := \widetilde{T}^{1/2}\left[\frac{G^w_{\theta_T,T}(\lambda)}{G^w_{\theta_T,T}(\pi)} - \frac{\lambda}{\pi}\right], \qquad \lambda \in [0, \pi],$$

Table 4

*Empirical size and power under local alternatives at 5% significance*

| | | $\boldsymbol{H_0}$: AR(1), $\boldsymbol{H_1}$: ARMA(1,1), $\boldsymbol{\eta = \tau/T^{1/2}}$ | | | | | | | | |
|---|---|---|---|---|---|---|---|---|---|---|
| $\tau$ | $\rho$ | $\hat{C}_T$ | $\hat{W}_{3,T}$ | $\hat{W}_{6,T}$ | $\lvert\hat{\psi}_{3,T}\rvert$ | $\lvert\hat{\psi}_{6,T}\rvert$ | $\hat{\psi}^+_{3,T}$ | $\hat{\psi}^+_{6,T}$ | $\tilde{Q}_{3,T}$ | $\tilde{Q}_{6,T}$ |
| | | | | | | $T = 200$ | | | | |
| 0 | 0.0 | 4.13 | 3.09 | 3.58 | 3.98 | 4.39 | 4.18 | 4.39 | 2.65 | 3.36 |
| | 0.5 | 3.62 | 2.80 | 2.22 | 3.68 | 4.04 | 3.93 | 4.14 | 2.67 | 3.31 |
| | 0.8 | 3.06 | 2.38 | 1.86 | 3.00 | 3.21 | 3.45 | 3.64 | 2.93 | 3.46 |
| 1 | 0.0 | 4.22 | 3.10 | 2.58 | 3.88 | 4.23 | 3.74 | 3.93 | 2.76 | 3.40 |
| | 0.5 | 5.52 | 4.08 | 2.90 | 5.51 | 5.76 | 8.86 | 9.20 | 3.08 | 3.61 |
| | 0.8 | 7.81 | 5.63 | 3.66 | 7.77 | 7.98 | 13.13 | 13.62 | 5.47 | 5.05 |
| 2 | 0.0 | 5.01 | 3.50 | 2.79 | 3.77 | 4.06 | 3.36 | 3.46 | 3.45 | 3.82 |
| | 0.5 | 8.53 | 6.10 | 4.02 | 8.58 | 9.06 | 14.33 | 14.61 | 4.51 | 4.56 |
| | 0.8 | 18.07 | 13.73 | 8.53 | 20.63 | 21.26 | 30.93 | 31.41 | 12.52 | 10.66 |
| 3 | 0.0 | 7.79 | 5.04 | 3.76 | 4.62 | 4.92 | 6.00 | 6.06 | 5.60 | 5.32 |
| | 0.5 | 10.64 | 7.80 | 5.16 | 10.84 | 11.25 | 17.39 | 17.87 | 5.76 | 5.41 |
| | 0.8 | 32.10 | 27.17 | 17.65 | 37.68 | 38.18 | 50.25 | 50.49 | 23.84 | 20.09 |
| 4 | 0.0 | 14.60 | 9.51 | 6.65 | 10.86 | 11.01 | 16.70 | 16.78 | 11.03 | 8.99 |
| | 0.5 | 10.67 | 8.16 | 5.42 | 10.65 | 11.01 | 17.11 | 17.57 | 5.93 | 5.56 |
| | 0.8 | 45.29 | 42.62 | 29.55 | 52.48 | 52.79 | 64.96 | 64.97 | 36.18 | 31.63 |
| | | | | | | $T = 500$ | | | | |
| 0 | 0.0 | 4.70 | 4.43 | 3.86 | 4.66 | 5.68 | 4.52 | 4.62 | 2.99 | 3.64 |
| | 0.5 | 4.50 | 4.23 | 3.70 | 4.53 | 4.55 | 4.50 | 4.52 | 2.99 | 3.80 |
| | 0.8 | 4.39 | 3.94 | 3.40 | 4.22 | 4.26 | 4.37 | 4.38 | 3.34 | 3.78 |
| 1 | 0.0 | 4.74 | 4.37 | 3.83 | 4.70 | 4.75 | 4.31 | 4.35 | 3.02 | 3.70 |
| | 0.5 | 6.68 | 5.72 | 4.73 | 6.71 | 6.61 | 10.25 | 10.36 | 3.75 | 4.32 |
| | 0.8 | 9.56 | 8.06 | 6.00 | 10.03 | 10.08 | 16.20 | 16.28 | 6.26 | 5.82 |
| 2 | 0.0 | 5.00 | 4.47 | 3.90 | 4.76 | 4.87 | 3.61 | 3.62 | 3.34 | 3.90 |
| | 0.5 | 11.06 | 8.94 | 6.81 | 11.48 | 11.43 | 18.23 | 18.17 | 6.06 | 5.88 |
| | 0.8 | 23.21 | 19.66 | 13.89 | 26.87 | 26.88 | 38.01 | 37.99 | 15.66 | 13.35 |
| 3 | 0.0 | 6.31 | 5.17 | 4.38 | 4.95 | 5.03 | 3.19 | 3.18 | 4.25 | 4.55 |
| | 0.5 | 16.44 | 13.17 | 9.58 | 17.26 | 17.24 | 26.26 | 26.03 | 9.45 | 8.39 |
| | 0.8 | 42.78 | 38.92 | 28.30 | 50.11 | 49.91 | 62.36 | 62.42 | 32.23 | 27.37 |
| 4 | 0.0 | 9.48 | 6.98 | 5.57 | 5.09 | 5.16 | 4.09 | 4.07 | 6.40 | 5.98 |
| | 0.5 | 21.08 | 17.22 | 12.42 | 22.10 | 21.95 | 32.15 | 31.99 | 12.84 | 10.89 |
| | 0.8 | 62.44 | 60.69 | 47.41 | 70.99 | 70.86 | 80.69 | 80.67 | 52.01 | 46.42 |

$\lvert\hat{\psi}_{n,T}\rvert$ denotes two-sided tests, whereas $\hat{\psi}^+_{n,T}$ are one-sided (right-hand side) tests.

where $G^w_{\theta,T}$ is analogous to $G_{\theta,T}$, but using the tapered periodogram, for example,

$$I^w_X(\lambda) := \frac{\lvert\sum_{t=1}^T w(t)X(t)e^{it\lambda}\rvert^2}{2\pi\sum_{t=1}^T w^2(t)}.$$

Here $\theta_T = \arg\min_{\theta\in\Theta} G^w_{\theta,T}(\pi)$ is the Whittle estimator proposed by Velasco and Robinson [45], which admits a similar asymptotic first order expansion as in (8), and where $w$ is a taper function, for example, the full cosine taper

$$w(t) = \frac{1}{2}\left(1 - \cos\left(\frac{2\pi t}{T}\right)\right), \qquad t = 1, \dots, T.$$

If the full cosine taper is used, because of its desirable asymptotic properties (see [44]), it is recommended in practice to base our tests on the empirical process $\beta^w_{\theta_T,T}$, where

$$\beta^w_{\theta,T}(\lambda_m) := \left(\frac{P_4^2}{\widetilde{T}}\right)^{1/2} \frac{2\pi}{G^w_{\theta,T}(\pi)} \sum_{j=1}^{m} e^w_{\theta,T}(j), \qquad m = 1, \dots, \overline{T},$$

with

$$e^w_{\theta,T}(j) := \frac{I^w_X(\lambda_j)}{h_\theta(\lambda_j)} - \gamma'_\theta(\lambda_j) b^w_{\theta,T}(j),$$

$$b^w_{\theta,T}(j) := A^{-1}_{\theta,T}(j)\frac{1}{\widetilde{T}} \sum_{k=j+1}^{\widetilde{T}} \gamma_\theta(\lambda_k)\frac{I^w_X(\lambda_k)}{h_\theta(\lambda_k)}$$

and

$$P_4^2 := \lim_{T\to\infty} \frac{T\sum_{t=1}^T w^4(t)}{(\sum_{t=1}^T w^2(t))^2} = \frac{35}{18}.$$

Under appropriate regularity conditions, it can be proved using tools in [44] and [45] that $\beta^w_{\theta_T,T} \Rightarrow B_\pi$.

Finally, the methodology can be extended to test the correlation structure of the innovations of regression models (e.g., distributed-lags models) using the martingale part of the $U_p$-process based on the residuals. When $\mathbb{E}(z(t)u(s)) = 0$ for all $t, s$, where $\{z(t)\}_{t=1}^T$ are the regressors and $\{u(t)\}_{t=1}^T$ the error term, the residual $U_p$-process is asymptotically equivalent to the $U_p$-process based on the true innovations, and there is no need to use tests based on the martingale part of the $U_p$-process. When $\mathbb{E}(z(t)u(t-s)) \neq 0$ for some $s > 0$, the first-order expansion of the residual $U_p$-process depends on the cross-spectrum of the innovations and regressors. However, it seems possible to apply the results in this paper to implement tests based on the (approximate) martingale part of this $U_p$-process with estimated parameters.

**6. Lemmas.** This section provides a series of lemmas which will be used in the proofs of the main results. Some of them can be of independent interest. Henceforth, $z^{(k)}$ denotes the $k$th element of a $p\times 1$ vector $z$ and $K$ a finite positive constant. Also, we shall abbreviate $g(\lambda_j)$ by $g_j$ for a generic function $g(\lambda)$.

LEMMA 1. *Let* $\zeta:(0,\pi]\to\mathbb{R}^p$ *be a function such that* $\|\zeta(\lambda)\|\le K|\log\lambda|^{\ell}$, $\ell\ge 1$, *and* $\|\partial\zeta(\lambda)/\partial\lambda\|\le K\lambda^{-1}|\log\lambda|^{\ell-1}$ *for all* $\lambda>0$. *Then, as* $T\to\infty$,

$$\sup_{\lambda\in[0,\pi]}\left\|\frac{1}{\widetilde{T}}\sum_{j=1}^{[\widetilde{T}\lambda/\pi]}\zeta_j-\frac{1}{\pi}\int_0^{\lambda}\zeta(x)\,dx\right\|\le K\frac{(\log\widetilde{T})^{\ell}}{\widetilde{T}}. \tag{16}$$

PROOF. The left-hand side of (16) is bounded by

$$\sup_{\lambda\in[0,\pi/\widetilde{T})}\left\|\frac{1}{\pi}\int_0^{\lambda}\zeta(x)\,dx\right\|+\sup_{\lambda\in[\pi/\widetilde{T},\pi]}\left\|\frac{1}{\widetilde{T}}\sum_{j=1}^{[\widetilde{T}\lambda/\pi]}\zeta_j-\frac{1}{\pi}\int_0^{\lambda}\zeta(x)\,dx\right\|. \tag{17}$$

The first term of (17) is bounded by

$$\frac{1}{\pi}\int_0^{\pi/\widetilde{T}}\|\zeta(x)\|\,dx\le K\int_0^{\pi/\widetilde{T}}|\log x|^{\ell}\,dx\le K\frac{(\log\widetilde{T})^{\ell}}{\widetilde{T}}.$$

Next, by the triangle inequality, the second term of (17) is bounded by

$$\begin{aligned}&\sup_{\lambda\in[\pi/\widetilde{T},\pi]}\left\|\frac{1}{\widetilde{T}}\zeta(\lambda)-\frac{1}{\pi}\int_0^{\pi/\widetilde{T}}\zeta(x)\,dx\right\|\\&\quad+\sup_{\lambda\in[\pi/\widetilde{T},\pi]}\frac{1}{\pi}\sum_{j=1}^{[\widetilde{T}\lambda/\pi]-1}\int_{j\pi/\widetilde{T}}^{(j+1)\pi/\widetilde{T}}\|\zeta_j-\zeta(x)\|\,dx.\end{aligned} \tag{18}$$

The first term of (18) is bounded by $K\widetilde{T}^{-1}(\log\widetilde{T})^{\ell}$ since $\|\zeta(x)\|\le K|\log x|^{\ell}$. Next, by the mean value theorem, the second term of (18) is bounded by

$$\begin{aligned}K\sum_{j=1}^{\widetilde{T}-1}\int_{j\pi/\widetilde{T}}^{(j+1)\pi/\widetilde{T}}\frac{1}{\lambda_j}\left|\frac{j\pi}{\widetilde{T}}-x\right||\log x|^{\ell-1}\,dx&\le K\sum_{j=1}^{\widetilde{T}-1}\frac{1}{j}\int_{j\pi/\widetilde{T}}^{(j+1)\pi/\widetilde{T}}|\log x|^{\ell-1}\,dx\\&\le\frac{K(\log\widetilde{T})^{\ell}}{\widetilde{T}}.\end{aligned}$$

□

The next lemma corresponds to Giraitis, Hidalgo and Robinson's [14] Lemma 4.4, which we state without proof for easy reference. For this purpose, let $u_j:=h_j^{-1/2}(2\pi T)^{-1/2}\sum_{t=1}^{T}X(t)e^{it\lambda_j}$, $v_j:=(2\pi T)^{-1/2}\sum_{t=1}^{T}\varepsilon(t)e^{it\lambda_j}$ and $R_{X\varepsilon}(\lambda)$ be the spectral coherency ([6], pages 256–257) between $X$ and $\varepsilon$. Also, herewith $\overline{c}$ will denote the conjugate of the complex number $c$.

LEMMA 2. *Assuming* A1 *and* A2, *then, as* $T\to\infty$, *the following relations hold uniformly over* $1\le j<k\le\widetilde{T}$:

$$\mathbb{E}(u_j\overline{v}_j)=R_{X\varepsilon,j}+O(j^{-1}\log(j+1));$$

$$\mathbb{E}(u_j v_j) = O(j^{-1}\log(j+1));$$
$$\max(|\mathbb{E}(u_k \overline{v}_j)|, |\mathbb{E}(u_k v_j)|) = O(j^{-1}\log(k));$$
$$\max(|\mathbb{E}(v_k \overline{u}_j)|, |\mathbb{E}(v_k u_j)|) = O(j^{-1}\log(k)).$$

The next lemma corresponds to the proof of expression (4.8) of [37], pages 1648–1651, using the orders of magnitude of the terms $a_1, a_2, b_1$ and $b_2$ in [37] and Lemma 3 there, but using our Lemma 2 instead of Robinson's [36] Theorems 1 and 2 when appropriate.

LEMMA 3. *Let $\zeta : [0,\pi] \to \mathbb{R}^p$ satisfy the same conditions on $\phi_{\theta_0}$ in* A3(a)–(c). *Then, assuming* A1 *and* A2, *as $T \to \infty$, for $1 \le r < s \le \widetilde{T}$, $h = 1, \dots, p$:*

$$\mathbb{E}\left|\sum_{j=r}^{s} \zeta_j^{(h)} v_j(\overline{u}_j - \overline{v}_j)\right|^2$$
$$\le K \log^2(T) \sum_{j=r}^{s}\left\{ j^{-1}\log(T) + \sum_{k=r}^{s}(j^{-2}\log^2(T) + j^{-1}k^{-1/2})\right\}.$$

LEMMA 4. *Let $\zeta : [0,\pi] \to \mathbb{R}^p$ satisfy the same conditions on $\phi_{\theta_0}$ in* A3(a)–(c) *and write*

$$\alpha_T^{\zeta}(\lambda) := \frac{1}{\widetilde{T}^{1/2}} \sum_{j=1}^{[\widetilde{T}\lambda/\pi]} \zeta_j \left(I_{\varepsilon,j} - \frac{\sigma^2}{2\pi}\right)$$

*and*

$$\tilde{\alpha}_T^{\zeta}(\lambda) := \frac{1}{\widetilde{T}^{1/2}} \sum_{j=1}^{[\widetilde{T}\lambda/\pi]} \zeta_j \left(\frac{I_{X,j}}{h_j} - \frac{\sigma^2}{2\pi}\right).$$

*Then, under the conditions of Theorem 1, for some $0 < \delta < 1/6$,*

$$\mathbb{E} \sup_{\lambda\in[0,\pi]} \|\tilde{\alpha}_T^{\zeta}(\lambda) - \alpha_T^{\zeta}(\lambda)\| = O(T^{-\delta}). \tag{19}$$

PROOF. It suffices to show that (19) holds for each element of the vector $\tilde{\alpha}_T^{\zeta}(\lambda) - \alpha_T^{\zeta}(\lambda)$. Then, by the triangle inequality the left-hand side of (19) is bounded by

$$\begin{aligned}&\mathbb{E} \sup_{\lambda\in[0,\pi]} \frac{1}{\widetilde{T}^{1/2}} \sum_{j=1}^{[\widetilde{T}\lambda/\pi]} |\zeta_j^{(k)}||u_j - v_j|^2 \\ &\quad + 2\mathbb{E} \sup_{\lambda\in[0,\pi]} \left| \frac{1}{\widetilde{T}^{1/2}} \sum_{j=1}^{[\widetilde{T}\lambda/\pi]} \zeta_j^{(k)} v_j(\overline{u}_j - \overline{v}_j)\right|.\end{aligned} \tag{20}$$

The first term of (20) is bounded by

$$\begin{aligned}\frac{1}{\widetilde{T}^{1/2}}\sum_{j=1}^{\widetilde{T}}|\zeta_j^{(k)}|\biggl\{\biggl(\mathbb{E}|u_j|^2-\frac{\sigma^2}{2\pi}\biggr)-\biggl(\mathbb{E}(u_j\overline{v}_j)-\frac{\sigma^2}{2\pi}\biggr)\\ -\biggl(\mathbb{E}(\overline{u}_jv_j)-\frac{\sigma^2}{2\pi}\biggr)+\biggl(\mathbb{E}|v_j|^2-\frac{\sigma^2}{2\pi}\biggr)\biggr\}\\ =O\biggl(\frac{\log T}{\widetilde{T}^{1/2}}\sum_{j=1}^{\widetilde{T}}\frac{\log(j+1)}{j}\biggr)=O(T^{-\delta}),\end{aligned}$$

by Lemma 2, because $\mathbb{E}|v_j|^2=(2\pi)^{-1}\sigma^2$, and by assumption $|\zeta_j^{(k)}|\le K\log T$.

Next, to show that the second term of (20) is $O(T^{-\delta})$, it suffices to show that

$$\mathbb{E}\max_{s=1,\dots,\widetilde{T}}\biggl|\frac{1}{\widetilde{T}^{1/2}}\sum_{j=1}^{s}\zeta_j^{(k)}v_j(\overline{u}_j-\overline{v}_j)\biggr|=O(T^{-\delta}). \tag{21}$$

By the triangle inequality the left-hand side of (21) is bounded by

$$\mathbb{E}\max_{s=1,\dots,[\widetilde{T}^\beta]}\biggl|\frac{1}{\widetilde{T}^{1/2}}\sum_{j=1}^{s}\zeta_j^{(k)}v_j(\overline{u}_j-\overline{v}_j)\biggr|+\mathbb{E}\biggl|\frac{1}{\widetilde{T}^{1/2}}\sum_{j=1}^{[\widetilde{T}^\beta]}\zeta_j^{(k)}v_j(\overline{u}_j-\overline{v}_j)\biggr| \tag{22}$$

$$+\,\mathbb{E}\max_{s=[\widetilde{T}^\beta]+1,\dots,\widetilde{T}}\biggl|\frac{1}{\widetilde{T}^{1/2}}\sum_{j=[\widetilde{T}^\beta]+1}^{s}\zeta_j^{(k)}v_j(\overline{u}_j-\overline{v}_j)\biggr|, \tag{23}$$

where $\frac{1}{3}<\beta<\frac{1}{2}$. Using the inequality

$$\Bigl(\sup_p|c_p|\Bigr)^2=\sup_p|c_p|^2\le\sum_p|c_p|^2, \tag{24}$$

by the Cauchy–Schwarz inequality the square of (22) is bounded by

$$\frac{4}{\widetilde{T}}\sum_{s=1}^{[\widetilde{T}^\beta]}\mathbb{E}\biggl|\sum_{j=1}^{s}\zeta_j^{(k)}v_j(\overline{u}_j-\overline{v}_j)\biggr|^2=O(\widetilde{T}^{2\beta-1}\log^4 T)=O(T^{-2\delta})$$

using Lemma 3.

To complete the proof, we need to show that $(23)=O(T^{-\delta})$. To that end, let $q=0,\dots,[\widetilde{T}^\varsigma]-1$ with $\frac{1}{3}<\varsigma<\beta$. By the triangle inequality (23) is bounded by

$$\mathbb{E}\frac{1}{\widetilde{T}^{1/2}}\max_{s=[\widetilde{T}^\beta]+1,\dots,\widetilde{T}}\biggl|\biggl\{\sum_{j=[\widetilde{T}^\beta]+1}^{s}-\sum_{j=[\widetilde{T}^\beta]+1}^{[\widetilde{T}^\beta]+q(s)\widetilde{T}/[\widetilde{T}^\varsigma]}\biggr\}\zeta_j^{(k)}v_j(\overline{u}_j-\overline{v}_j)\biggr|$$

(25)
$$+\mathbb{E}\frac{1}{\widetilde{T}^{1/2}}\max_{s=[\widetilde{T}^{\beta}]+1,\dots,\widetilde{T}}\left|\sum_{j=[\widetilde{T}^{\beta}]+1}^{[\widetilde{T}^{\beta}]+q(s)\widetilde{T}/[\widetilde{T}^{\varsigma}]}\zeta_j^{(k)}v_j(\overline{u}_j-\overline{v}_j)\right|,$$

where $q(s)$ denotes the value of $q=0,\dots,[\widetilde{T}^{\varsigma}]-1$ such that $[\widetilde{T}^{\beta}]+q(s)\widetilde{T}/[\widetilde{T}^{\varsigma}]$ is the largest integer smaller than or equal to $s$, and using the convention $\sum_c^d\equiv 0$ if $d<c$.

By the definition of $q(s)$ and the Cauchy–Schwarz inequality, the square of the second term of (25) is bounded by

$$\mathbb{E}\frac{1}{\widetilde{T}}\max_{q=0,\dots,[\widetilde{T}^{\varsigma}]-1}\left|\sum_{j=[\widetilde{T}^{\beta}]+1}^{[\widetilde{T}^{\beta}]+q\widetilde{T}/[\widetilde{T}^{\varsigma}]}\zeta_j^{(k)}v_j(\overline{u}_j-\overline{v}_j)\right|^2$$
$$\le\frac{1}{\widetilde{T}}\sum_{q=0}^{[\widetilde{T}^{\varsigma}]-1}\mathbb{E}\left|\sum_{j=[\widetilde{T}^{\beta}]+1}^{[\widetilde{T}^{\beta}]+q\widetilde{T}/[\widetilde{T}^{\varsigma}]}\zeta_j^{(k)}v_j(\overline{u}_j-\overline{v}_j)\right|^2$$

by (24). But, using Lemma 3, we have that the right-hand side of the last displayed inequality is bounded by

$$K\frac{\log^4 T}{\widetilde{T}}\sum_{q=0}^{[\widetilde{T}^{\varsigma}]-1}\left(1+\frac{|q|_+\widetilde{T}^{1-\varsigma}}{\widetilde{T}^{\beta}}+|q|_+^{1/2}\widetilde{T}^{1/2(1-\varsigma)}\right)$$
$$\le K\log^4 T(\widetilde{T}^{\varsigma-\beta}+\widetilde{T}^{\varsigma-1/2})\le K\widetilde{T}^{-2\delta},$$

where $|q|_+=\max\{1,|q|\}$. To complete the proof, we need to show that the first term in (25) is $O(T^{-\delta})$. To that end, we note that this term is bounded by

$$\mathbb{E}\frac{1}{\widetilde{T}^{1/2}}\max_{q=0,\dots,[\widetilde{T}^{\varsigma}]-1}\max_{s}\left|\sum_{j=1+[\widetilde{T}^{\beta}]+q\widetilde{T}/[\widetilde{T}^{\varsigma}]}^{s}\zeta_j^{(k)}v_j(\overline{u}_j-\overline{v}_j)\right|,$$

where the $\max_s$ runs for all values $s=1+[\widetilde{T}^{\beta}]+q\widetilde{T}/[\widetilde{T}^{\varsigma}],\dots,[\widetilde{T}^{\beta}]+(q+1)\widetilde{T}/[\widetilde{T}^{\varsigma}]$. By the Cauchy–Schwarz inequality and (24), the square of the last displayed expression is bounded by

$$\frac{1}{\widetilde{T}}\sum_{q=0}^{[\widetilde{T}^{\varsigma}]-1}\sum_{s=1+[\widetilde{T}^{\beta}]+q\widetilde{T}/[\widetilde{T}^{\varsigma}]}^{[\widetilde{T}^{\beta}]+(q+1)\widetilde{T}/[\widetilde{T}^{\varsigma}]}\mathbb{E}\left|\sum_{j=1+[\widetilde{T}^{\beta}]+q\widetilde{T}/[\widetilde{T}^{\varsigma}]}^{s}\zeta_j^{(k)}v_j(\overline{u}_j-\overline{v}_j)\right|^2$$
$$\le K\frac{\log^4\widetilde{T}}{\widetilde{T}}\sum_{q=0}^{[\widetilde{T}^{\varsigma}]-1}\sum_{s=1+[\widetilde{T}^{\beta}]+q\widetilde{T}/[\widetilde{T}^{\varsigma}]}^{[\widetilde{T}^{\beta}]+(q+1)\widetilde{T}/[\widetilde{T}^{\varsigma}]}\left\{\frac{1}{|q|_+}+\frac{\widetilde{T}^{(1-\varsigma)/2}}{|q|_+^{3/2}}\right\}$$

$$\leq K\frac{\log^4 \widetilde{T}}{\widetilde{T}}(\widetilde{T}^{1-\varsigma}\log T+\widetilde{T}^{3(1-\varsigma)/2})\leq K\widetilde{T}^{(1-3\varsigma)/2}\log^4 T\leq K\widetilde{T}^{-2\delta},$$

where in the first inequality we have used Lemma 3 and that, for $q\geq 1$ and $\psi\geq 0$,

$$\sum_{j=1+[\widetilde{T}^{\beta}]+q\widetilde{T}/[\widetilde{T}^{\varsigma}]}^{s} j^{-\psi}\leq \frac{K}{(\widetilde{T}^{\beta}+q\widetilde{T}^{1-\varsigma})^{\psi}}\left(\sum_{j=1+[\widetilde{T}^{\beta}]+q\widetilde{T}/[\widetilde{T}^{\varsigma}]}^{[\widetilde{T}^{\beta}]+(q+1)\widetilde{T}/[\widetilde{T}^{\varsigma}]} 1\right)$$
$$\leq \frac{K\widetilde{T}^{(1-\varsigma)(1-\psi)}}{q^{\psi}}.$$

This completes the proof. $\square$

REMARK 1. Lemma 4 holds for $\alpha_T^{\zeta}(\lambda)$ and $\tilde{\alpha}_T^{\zeta}(\lambda)$ replaced by

$$\ddot{\alpha}_T^{\zeta}(\lambda):=\alpha_T^{\zeta}(\pi)-\alpha_T^{\zeta}(\lambda),\qquad \ddot{\tilde{\alpha}}_T^{\zeta}(\lambda):=\tilde{\alpha}_T^{\zeta}(\pi)-\tilde{\alpha}_T^{\zeta}(\lambda),$$

respectively. This is so because the triangle inequality implies that

$$\mathbb{E}\sup_{\lambda\in[0,\pi]}|\ddot{\alpha}_T^{\zeta}(\lambda)-\ddot{\tilde{\alpha}}_T^{\zeta}(\lambda)|\leq 2\mathbb{E}\sup_{\lambda\in[0,\pi]}|\alpha_T^{\zeta}(\lambda)-\tilde{\alpha}_T^{\zeta}(\lambda)|.$$

Define, for $\mu$ and $\vartheta\in[0,\pi]$,

$$c_s(\mu,\vartheta)=\frac{2}{T\widetilde{T}^{1/2}}\sum_{p=[\widetilde{T}\mu/\pi]+1}^{[\widetilde{T}\vartheta/\pi]}\zeta_p\cos(s\lambda_p), \tag{26}$$

where $\zeta$ is as in Lemma 1 and $\mu<\vartheta$.

LEMMA 5. *For* $0\leq\mu<\vartheta_1,\vartheta_2\leq\pi$, *as* $T\to\infty$,

$$\sum_{t=1}^{T-1}\sum_{s=1}^{T-t}c_s(\mu,\vartheta_1)c_s'(\mu,\vartheta_2)=g(\mu,\vartheta_1,\vartheta_2)(1+o(1)), \tag{27}$$

*where* $g(\mu,\vartheta_1,\vartheta_2)=\pi^{-1}\int_{\mu}^{\vartheta_1\wedge\vartheta_2}\zeta(u)\zeta'(u)\,du-(\pi^{-1}\int_{\mu}^{\vartheta_1}\zeta(u)\,du)(\pi^{-1}\times\int_{\mu}^{\vartheta_2}\zeta'(u)\,du)$.

PROOF. A typical component of the matrix on the left-hand side of (27) is

$$\frac{4}{T^2\widetilde{T}}\sum_{p_1=[\widetilde{T}\mu/\pi]+1}^{[\widetilde{T}\vartheta_1/\pi]}\zeta_{p_1}^{(k_1)}\sum_{p_2=[\widetilde{T}\mu/\pi]+1}^{[\widetilde{T}\vartheta_2/\pi]}\zeta_{p_2}^{(k_2)}\sum_{t=1}^{T-1}\sum_{s=1}^{T-t}\cos(s\lambda_{p_1})\cos(s\lambda_{p_2})$$

$$
\begin{aligned}
&= \frac{4}{T^2\widetilde{T}} \sum_{p=[\widetilde{T}\mu/\pi]+1}^{[\widetilde{T}\vartheta_1/\pi]\wedge[\widetilde{T}\vartheta_2/\pi]} \zeta_p^{(k_1)}\zeta_p^{(k_2)} \sum_{t=1}^{T-1}\sum_{s=1}^{T-t}\cos^2(s\lambda_p) \\
&\quad + \frac{2}{T^2\widetilde{T}} \sum_{p_1=[\widetilde{T}\mu/\pi]+1}^{[\widetilde{T}\vartheta_1/\pi]} \zeta_{p_1}^{(k_1)} \\
&\quad \times \sum_{\substack{p_2=[\widetilde{T}\mu/\pi]+1\\ p_2\neq p_1}}^{[\widetilde{T}\vartheta_2/\pi]} \zeta_{p_2}^{(k_2)} \sum_{t=1}^{T-1}\sum_{s=1}^{T-t}\{\cos(s\lambda_{p_1+p_2})+\cos(s\lambda_{p_1-p_2})\}.
\end{aligned} \tag{28}
$$

Because $\cos^2\lambda = (1+\cos(2\lambda))/2$, then using formulae in [6], page 13, we have that $\sum_{t=1}^{T-1}\sum_{s=1}^{T-t}\cos^2(s\lambda_p) = (T-1)^2/4$ and, for $p_1 \neq p_2$,

$$
\sum_{t=1}^{T-1}\sum_{s=1}^{T-t}\{\cos(s\lambda_{p_1+p_2})+\cos(s\lambda_{p_1-p_2})\} = -T
$$

and, hence, we conclude that the right-hand side of (28) is, recalling that $\widetilde{T} = [T/2]$,

$$
\begin{aligned}
&\frac{(T-1)^2}{T^2}\left(\frac{1}{\widetilde{T}} \sum_{p=[\widetilde{T}\mu/\pi]+1}^{[\widetilde{T}\vartheta_1/\pi]\wedge[\widetilde{T}\vartheta_2/\pi]} \zeta_p^{(k_1)}\zeta_p^{(k_2)}\right) \\
&\quad - \frac{2}{T\widetilde{T}} \sum_{p_1=[\widetilde{T}\mu/\pi]+1}^{[\widetilde{T}\vartheta_1/\pi]} \zeta_{p_1}^{(k_1)} \sum_{\substack{p_2=[\widetilde{T}\mu/\pi]+1\\ p_2\neq p_1}}^{[\widetilde{T}\vartheta_2/\pi]} \zeta_{p_2}^{(k_2)} \\
&= g^{(k_1,k_2)}(\mu,\vartheta_1,\vartheta_2)(1+o(1)),
\end{aligned}
$$

by Lemma 1 and where $g^{(k_1,k_2)}(\mu,\vartheta_1,\vartheta_2)$ denotes the $(k_1,k_2)$th element of the matrix $g(\mu,\vartheta_1,\vartheta_2)$. $\square$

We now introduce the following notation. For $0 \le v_1 < v_2 \le \pi$,

$$
\mathcal{E}_{1,T}(v_1,v_2) := \left(\frac{1}{\widetilde{T}} \sum_{p=[\widetilde{T}v_1/\pi]+1}^{[\widetilde{T}v_2/\pi]} \zeta_p\right)\left(\frac{\widetilde{T}^{1/2}}{T}\sum_{t=1}^{T}(\varepsilon^2(t)-\sigma^2)\right), \tag{29}
$$

$$
\mathcal{E}_{2,T}(v_1,v_2) := \sum_{t=2}^{T}\varepsilon(t)\sum_{s=1}^{t-1}\varepsilon(s)c_{t-s}(v_1,v_2), \tag{30}
$$

where $c_t(\cdot,\cdot)$ is given in (26) and $\zeta$ is as in Lemma 1.

LEMMA 6. *Let $0 \le v_1 < v < v_2 < \pi$. Then assuming* A1, *for $k = 1, \dots, p$ and for some $\beta > 0$ and $0 < \delta < 1$,*

$$\mathbb{E}(|\mathcal{E}_{j,T}^{(k)}(v_1, v)|^\beta |\mathcal{E}_{j,T}^{(k)}(v, v_2)|^\beta) \le K(v_2 - v_1)^{2-\delta}, \qquad j = 1, 2, \tag{31}$$

*where $\mathcal{E}_{1,T}^{(k)}(v_1, v)$ and $\mathcal{E}_{2,T}^{(k)}(v_1, v)$ are the kth components of* (29) *and* (30), *respectively.*

PROOF. We begin with $j = 1$. By Lemma 1,

$$\left| \frac{1}{\widetilde{T}} \sum_{p=[\widetilde{T}v_1/\pi]+1}^{[\widetilde{T}v_2/\pi]} \zeta_p^{(k)} - \frac{1}{\pi} \int_{v_1}^{v_2} \zeta^{(k)}(x)\,dx \right| \le K \frac{|\log \widetilde{T}|^\ell}{\widetilde{T}} \le K(v_2 - v_1)^{1-\delta/2},$$

after we notice that we can take $\widetilde{T}^{-1} \le (v_2 - v_1)$, since otherwise (31) holds trivially. On the other hand, A1 implies that $\mathbb{E}(\sum_{t=1}^T (\varepsilon^2(t) - \sigma^2))^2 \le KT$. So, using the inequality $(v_2 - v)(v - v_1) < (v_2 - v_1)^2$ and the Cauchy–Schwarz inequality, we have that $\mathbb{E}(|\mathcal{E}_{1,T}^{(k)}(v_1, v)||\mathcal{E}_{1,T}^{(k)}(v, v_2)|) \le K(v_2 - v_1)^{2-\delta}$.

To complete the proof, it suffices to examine that the inequality in (31) holds for $j = 2$. Now

$$\mathbb{E}(\mathcal{E}_{2,T}^{(k)}(v_1, v_2))^4 = 16 \prod_{j=1}^4 \sum_{1 \le s_j < t_j \le T} c_{t_j - s_j}^{(k)}(v_1, v_2) \mathbb{E}(\varepsilon(t_1)\varepsilon(s_1) \dots \varepsilon(t_4)\varepsilon(s_4)).$$

Since the number of equal indices in the set $\{t_1, s_1, \dots, t_4, s_4\}$ does not exceed 4, by assumption A1 it follows that $|\mathbb{E}(\varepsilon(t_1)\varepsilon(s_1) \dots \varepsilon(t_4)\varepsilon(s_4))| \le K$. Moreover, by A1 the inequality $|\mathbb{E}(\varepsilon(t_1)\varepsilon(s_1) \dots \varepsilon(t_4)\varepsilon(s_4))| \ne 0$ can hold only if any $t_j, s_j$ are repeated in $\{t_1, s_1, \dots, t_4, s_4\}$ at least twice. Hence, by the Cauchy–Schwarz inequality, we obtain that

$$\begin{aligned} \mathbb{E}(\mathcal{E}_{2,T}^{(k)}(v_1, v_2))^4 &\le K \prod_{j=1}^4 \left( \sum_{1 \le s_j < t_j \le T} (c_{t_j - s_j}^{(k)}(v_1, v_2))^2 \right)^{1/2} \\ &= K \left( \sum_{1 \le s < t \le T} (c_{t-s}^{(k)}(v_1, v_2))^2 \right)^2. \end{aligned}$$

But by Lemma 5 the right-hand side of the last displayed equation is bounded by

$$K \left( \frac{1}{\pi} \int_{v_1}^{v_2} (\zeta^{(k)}(u))^2\,du - \left( \frac{1}{\pi} \int_{v_1}^{v_2} \zeta^{(k)}(u)\,du \right)^2 \right)^2 \le K(v_2 - v_1)^{2-\delta}$$

because $|\int_{v_1}^{v_2} (\zeta^{(k)}(x))^p\,dx| \le K|v_2 - v_1|^{1-\delta/2}$ for $p = 1, 2$. This concludes the proof choosing $\beta = 2$ by the Cauchy–Schwarz inequality. $\square$

LEMMA 7. *Denote* $\eta_p := I_{\varepsilon,p} - \sigma^2/(2\pi)$ *and*

$$R_T^1(v) = \frac{2\pi}{\widetilde{T}^{1/2}} \sum_{p=1}^{[\widetilde{T}v/\pi]} \zeta_p \eta_p \quad \text{and} \quad R_T^2(v) = \frac{2\pi}{\widetilde{T}^{1/2}} \sum_{p=[\widetilde{T}v/\pi]+1}^{\widetilde{T}} \zeta_p \eta_p, \qquad (0 \le v < \pi)$$

*with* $\zeta$ *as in Lemma* 1. *Let* $0 \le v_1 < v < v_2 \le \pi$. *Then assuming* A1, *for some* $\beta > 0$ *and* $0 < \delta < 1$:

$$\text{(a)} \qquad \mathbb{E}(\|R_T^j(v_2) - R_T^j(v)\|^\beta \|R_T^j(v) - R_T^j(v_1)\|^\beta) \le K(v_2 - v_1)^{2-\delta}, \qquad j = 1, 2. \tag{32}$$

(b) $R_T^j(v) \xrightarrow{d} \mathcal{N}(0, 4\pi^2 V^{(j)}(v)), \qquad j = 1, 2,$

*where* $V^{(1)}(v) = \sigma^4 \int_0^v \zeta(u)\zeta'(u)\,du/\pi + \sigma^4\kappa \int_0^v \zeta(u)\,du \int_0^v \zeta'(u)\,du/\pi^2$ *and* $V^{(2)}(v) = \sigma^4 \int_v^\pi \zeta(u)\zeta'(u)\,du/\pi + \sigma^4\kappa \int_v^\pi \zeta(u)\,du \int_v^\pi \zeta'(u)\,du/\pi^2$, *with* $\kappa$ *denoting the fourth cumulant of* $\{\varepsilon(t)/\sigma\}_{t\in\mathbb{Z}}$.

PROOF. We begin with (a). We shall consider $R_T^2(v)$ only, $R_T^1(v)$ being similarly handled. From the definition of $\eta_p$, and

$$R_T^2(v) - R_T^2(v_2) = \frac{2\pi}{\widetilde{T}^{1/2}} \sum_{p=[\widetilde{T}v/\pi]+1}^{[\widetilde{T}v_2/\pi]} \zeta_p \eta_p,$$

we have that

$$R_T^2(v) - R_T^2(v_2) = \mathcal{E}_{1,T}(v, v_2) + \mathcal{E}_{2,T}(v, v_2),$$

where $\mathcal{E}_{1,T}(v, v_2)$ and $\mathcal{E}_{2,T}(v, v_2)$ are given in (29) and (30), respectively. Now (32) follows immediately from Lemma 6 and standard inequalities.

Part (b). We will examine $R_T^1(v) \xrightarrow{d} \mathcal{N}(0, 4\pi^2 V^{(1)}(v))$, the proof for $j = 2$ being handled identically. But this follows by an obvious extension of Theorem 4.2 of [14] because $\zeta(u)$ satisfies the same conditions on $h_n(u)$ there. □

LEMMA 8. *Assume* A1–A4. *Then we have that, for some* $0 < \delta < 1/6$,

$$\text{(a)} \quad \frac{2\pi}{\widetilde{T}^{1/2}} \sum_{j=1}^{[\widetilde{T}\lambda/\pi]} \zeta_j \left( \frac{I_{X,j}}{h_{\theta_T,j}} - \frac{\sigma^2}{2\pi} \right) = \frac{2\pi}{\widetilde{T}^{1/2}} \sum_{j=1}^{[\widetilde{T}\lambda/\pi]} \zeta_j \left( I_{\varepsilon,j} - \frac{\sigma^2}{2\pi} \right) - \left( \frac{\sigma^2}{\widetilde{T}} \sum_{j=1}^{[\widetilde{T}\lambda/\pi]} \zeta_j \phi'_{\theta_0,j} \right) \widetilde{T}^{1/2}(\theta_T - \theta_0) \tag{33}$$

$$+ O_p\left(\frac{1}{T^\delta}\right),$$

(b) $$\frac{2\pi}{\widetilde{T}^{1/2}} \sum_{j=[\widetilde{T}\lambda/\pi]+1}^{\widetilde{T}} \zeta_j \left(\frac{I_{X,j}}{h_{\theta_T,j}} - \frac{\sigma^2}{2\pi}\right)$$
$$= \frac{2\pi}{\widetilde{T}^{1/2}} \sum_{j=[\widetilde{T}\lambda/\pi]+1}^{\widetilde{T}} \zeta_j \left(I_{\varepsilon,j} - \frac{\sigma^2}{2\pi}\right)$$
$$- \left(\frac{\sigma^2}{\widetilde{T}} \sum_{j=[\widetilde{T}\lambda/\pi]+1}^{\widetilde{T}} \zeta_j \phi'_{\theta_0,j}\right) \widetilde{T}^{1/2}(\theta_T - \theta_0) + O_p\left(\frac{1}{T^\delta}\right),$$

*where the $O_p(1/T^\delta)$ terms are uniform in $\lambda \in [0,\pi]$, and where $\zeta(u)$ and $\|\zeta(u)\|$ are as in Lemma 1.*

PROOF. We examine (a), part (b) being handled similarly. The difference between the left-hand side of (33) and the first term on its right-hand side is

$$\frac{2\pi}{\widetilde{T}^{1/2}} \sum_{j=1}^{[\widetilde{T}\lambda/\pi]} \zeta_j \frac{I_{X,j}}{h_{\theta_0,j}} \left[\frac{h_{\theta_0,j}}{h_{\theta_T,j}} - 1 + \phi'_{\theta_0,j}(\theta_T - \theta_0)\right]$$
(34)
$$+ \frac{2\pi}{\widetilde{T}^{1/2}} \sum_{j=1}^{[\widetilde{T}\lambda/\pi]} \zeta_j \left(\frac{I_{X,j}}{h_{\theta_0,j}} - I_{\varepsilon,j}\right) - \frac{2\pi}{\widetilde{T}^{1/2}} \sum_{j=1}^{[\widetilde{T}\lambda/\pi]} \zeta_j \phi'_{\theta_0,j} \frac{I_{X,j}}{h_{\theta_0,j}} (\theta_T - \theta_0).$$

First we notice that

(35) $$\theta_T - \theta_0 = O_p(T^{-1/2}),$$

which follows by (8) in assumption A4, and because

(36) $$\frac{1}{\widetilde{T}^{1/2}} \sum_{k=1}^{\widetilde{T}} \phi_{\theta_0,k} \left(\frac{I_{X,k}}{h_{\theta_0,k}} - I_{\varepsilon,k}\right) = O_p(T^{-\delta})$$

(recall that under $H_0$, $h_j = h_{\theta_0,j}$), by Lemma 4 and Markov's inequality, and

$$\frac{2\pi}{\sigma^2 \widetilde{T}^{1/2}} \sum_{k=1}^{\widetilde{T}} \phi_{\theta_0,k} I_{\varepsilon,k} \xrightarrow{d} \mathcal{N}\left(0, \frac{1}{\pi} \int_0^\pi \phi_{\theta_0}(u) \phi'_{\theta_0}(u)\, du\right)$$
(37)
$$\stackrel{d}{=} \int_0^\pi \phi_{\theta_0}(u) B_\pi(du)$$

by Lemma 7 with $\zeta(u) = \phi_{\theta_0}(u)$. Notice also that $\sum_{k=1}^{\widetilde{T}} \phi_{\theta_0,k} = O(\log T)$ by Lemma 1 because (9) and A3 part (c) implies that $\phi_{\theta_0}(\lambda)$ satisfies the same conditions on $\zeta(\lambda)$ in Lemma 1.

Next, A3 part (d) implies that, uniformly in $\lambda \in [0,\pi]$, the norm of the first term of (34) is bounded by

$$K\widetilde{T}^{1/2}\|\theta_T - \theta_0\|^2 \frac{1}{\widetilde{T}} \sum_{j=1}^{[\widetilde{T}\lambda/\pi]} |\log^2 \lambda_j| \|\zeta_j\| \frac{I_{X,j}}{h_{\theta_0,j}} = O_p(T^{-1/2}), \tag{38}$$

because (35) implies that we can take $\delta = KT^{-1/2}$ in A3 part (d) so that $\lambda_j^{-\delta} < K$ when $\delta < KT^{-1/2}$ and $j \geq 1$, and also because by Markov's inequality and Lemmas 4 and 7,

$$\sup_{\lambda \in [0,\pi]} \left| \frac{1}{\widetilde{T}} \sum_{j=1}^{[\widetilde{T}\lambda/\pi]} |\log^2 \lambda_j| \|\zeta_j\| \left( \frac{I_{X,j}}{h_{\theta_0,j}} - \frac{\sigma^2}{2\pi} \right) \right| = O_p(T^{-1/2}),$$

and because by Lemma 1 with $\|\zeta(u)\| |\log^2(u)|$ there,

$$\sup_{\lambda \in [0,\pi]} \left| \frac{1}{\widetilde{T}} \sum_{j=1}^{[\widetilde{T}\lambda/\pi]} |\log^2 \lambda_j| \|\zeta_j\| - \frac{1}{\pi} \int_0^\lambda |\log^2(u)| \|\zeta(u)\| \, du \right| = o(\widetilde{T}^{-1/2}).$$

The second term of (34) is $O_p(T^{-\delta})$ by Lemma 4 and Markov's inequality. Next, proceeding similarly as in (38), since $\zeta(\lambda)\phi'_{\theta_0}(\lambda)$ satisfies the same conditions as $\zeta(\lambda)|\log \lambda|$, the third term of (34) is $\widetilde{T}^{-1}\sigma^2 \sum_{j=1}^{[\widetilde{T}\lambda/\pi]} \zeta_j \phi'_{\theta_0,j} \widetilde{T}^{1/2}(\theta_T - \theta_0) + O_p(T^{-\delta})$, which concludes the proof. □

Lemma 9. *Assuming* A1, *for any* $0 \leq \upsilon < (1-\delta)/4$, *with* $\delta$ *as in Lemma* 7, *we have that, for all* $k = 1, \ldots, p$,

$$\text{(a)} \quad \mathbb{E}\left( \frac{\mathcal{E}_{1,T}^{(k)}(\lambda_1, \pi)}{(\pi - \lambda_1)^\upsilon} - \frac{\mathcal{E}_{1,T}^{(k)}(\lambda_2, \pi)}{(\pi - \lambda_2)^\upsilon} \right)^2 \leq K(\lambda_2 - \lambda_1)^{2-\delta-2\upsilon}, \tag{39}$$

$$\text{(b)} \quad \mathbb{E}\left( \frac{\mathcal{E}_{2,T}^{(k)}(\lambda_1, \pi)}{(\pi - \lambda_1)^\upsilon} - \frac{\mathcal{E}_{2,T}^{(k)}(\lambda_2, \pi)}{(\pi - \lambda_2)^\upsilon} \right)^4 \leq K(\lambda_2 - \lambda_1)^{2-\delta-4\upsilon} \tag{40}$$

*for all* $0 < \lambda_1 < \lambda_2 < \pi$, *and where* $\mathcal{E}_{1,T}^{(k)}(\lambda_1, \lambda_2)$ *and* $\mathcal{E}_{2,T}^{(k)}(\lambda_1, \lambda_2)$ *are given in* (29) *and* (30), *respectively.*

Proof. We begin with (b). By standard inequalities the left-hand side of (40) is bounded by

$$K\mathbb{E}\left( \frac{1}{(\pi - \lambda_1)^\upsilon} \mathcal{E}_{2,T}^{(k)}(\lambda_1, \lambda_2) \right)^4 + K\left( \frac{1}{(\pi - \lambda_1)^\upsilon} - \frac{1}{(\pi - \lambda_2)^\upsilon} \right)^4 \mathbb{E}(\mathcal{E}_{2,T}^{(k)}(\lambda_2, \pi))^4.$$

By Lemma 6, for any $0<\delta<1$, we have that the last displayed expression is bounded by

$$K\frac{(\lambda_2-\lambda_1)^{2-\delta}}{(\pi-\lambda_1)^{4\upsilon}}+K\left(\frac{1}{(\pi-\lambda_1)^{\upsilon}}-\frac{1}{(\pi-\lambda_2)^{\upsilon}}\right)^4(\pi-\lambda_2)^{2-\delta}. \tag{41}$$

Consider the case $\lambda_2-\lambda_1\le 2^{-1}(\pi-\lambda_2)$ first. By the mean value theorem (41) is

$$\begin{aligned}&K\frac{(\lambda_2-\lambda_1)^{2-\delta}}{(\pi-\lambda_1)^{4\upsilon}}\\&\quad+\frac{K}{(\pi-\lambda_1)^{4\upsilon}(\pi-\lambda_2)^{\delta+4\upsilon-2}}\frac{\upsilon^4(\lambda_2-\lambda_1)^4}{(\beta(\pi-\lambda_1)+(1-\beta)(\pi-\lambda_2))^{4-4\upsilon}}\\&\qquad\le K(\lambda_2-\lambda_1)^{2-\delta-4\upsilon}+K(\pi-\lambda_2)^{-\delta-4\upsilon-2}(\lambda_2-\lambda_1)^4,\end{aligned}$$

where $\beta=\beta(\lambda_1,\lambda_2)\in(0,1)$, and then because $\pi-\lambda_1>\lambda_2-\lambda_1$ and $\pi-\lambda_1\ge\pi-\lambda_2>0$. But the right-hand side of the last displayed inequality is bounded by $K(\lambda_2-\lambda_1)^{2-\delta-4\upsilon}$ since $\lambda_2-\lambda_1\le 2^{-1}(\pi-\lambda_2)$.

Next, consider the case for which $2^{-1}(\pi-\lambda_2)<\lambda_2-\lambda_1$. Using the inequality $a^{\varsigma}-b^{\varsigma}\le(a-b)^{\varsigma}$ for any $0<\varsigma<1$ and $a\ge b$, we have that (41) is bounded by

$$K(\lambda_2-\lambda_1)^{2-\delta-4\upsilon}+K\frac{(\lambda_2-\lambda_1)^{4\upsilon}(\pi-\lambda_2)^{2-\delta}}{(\pi-\lambda_1)^{4\upsilon}(\pi-\lambda_2)^{4\upsilon}}\le K(\lambda_2-\lambda_1)^{2-\delta-4\upsilon},$$

where we have used $0<\lambda_2-\lambda_1\le\pi-\lambda_1$ and $\pi-\lambda_2<2(\lambda_2-\lambda_1)$. This completes the proof of part (b).

Next part (a). By definition and A1, the left-hand side of (39) is bounded by

$$\begin{aligned}&\frac{K}{(\pi-\lambda_1)^{2\upsilon}}\left(\frac{1}{\widetilde{T}}\sum_{j=[\widetilde{T}\lambda_1/\pi]+1}^{[\widetilde{T}\lambda_2/\pi]}\zeta_j^{(k)}\right)^2\\&\quad+K\left(\frac{1}{(\pi-\lambda_1)^{\upsilon}}-\frac{1}{(\pi-\lambda_2)^{\upsilon}}\right)^2\left(\frac{1}{\widetilde{T}}\sum_{j=[\widetilde{T}\lambda_2/\pi]+1}^{\widetilde{T}}\zeta_j^{(k)}\right)^2\\&\qquad\le K(\lambda_2-\lambda_1)^{2-\delta-2\upsilon}\end{aligned}$$

by Lemma 1, and then proceed as in part (b). □

In what follows we shall abbreviate $\gamma'_{\theta,q}A^{-1}_{\theta,T}(q)$ by $H_{\theta,T}(q)$.

LEMMA 10. *Assuming* A1–A5, *for all* $\epsilon > 0$,

$$\lim_{\lambda_0 \to \pi} \limsup_{T \to \infty} \Pr\left\{ \sup_{\lambda_0 \le \lambda \le \pi} \left| \frac{1}{\widetilde{\widetilde{T}}} \sum_{k=[\overline{T}\lambda_0/\pi]+1}^{[\overline{T}\lambda/\pi]} \frac{H_{\theta_0,T}(k)}{\widetilde{T}^{1/2}} \right.\right. \tag{42}$$
$$\left.\left. \times \sum_{j=k+1}^{\widetilde{T}} \gamma_{\theta_0,j} \left( \frac{I_{X,j}}{h_{\theta_T,j}} - \frac{\sigma^2}{2\pi} \right) \right| > \epsilon \right\} = 0.$$

PROOF. Abbreviate $h_{\theta_T,j}^{-1} I_{X,j} - I_{\varepsilon,j}$ by $\varkappa_j$ and take $\lambda_0 > \pi/2$ without loss of generality. Noting that $h_{\theta_T,j}^{-1} I_{X,j} - \sigma^2/(2\pi) = \varkappa_j + \eta_j$, where $\eta_j = I_{\varepsilon,j} - \sigma^2/(2\pi)$, we have

$$\sup_{\lambda_0 \le \lambda \le \pi} \left| \frac{1}{\widetilde{\widetilde{T}}} \sum_{k=[\overline{T}\lambda_0/\pi]+1}^{[\overline{T}\lambda/\pi]} \frac{H_{\theta_0,T}(k)}{\widetilde{T}^{1/2}} \sum_{j=k+1}^{\widetilde{T}} \gamma_{\theta_0,j} (\varkappa_j + \eta_j) \right|$$
$$\le \frac{K}{\widetilde{\widetilde{T}}} \sum_{k=[\overline{T}\lambda_0/\pi]+1}^{\overline{T}} \|H_{\theta_0,T}(k)\| \left(1 - \frac{k}{\widetilde{T}}\right)^{\delta/2} \tag{43}$$
$$\times \left\{ \sup_{[\overline{T}\lambda_0/\pi] \le k \le \overline{T}} \left\| \frac{(1 - k/\widetilde{T})^{-\delta/2}}{\widetilde{T}^{1/2}} \sum_{j=k+1}^{\widetilde{T}} \gamma_{\theta_0,j} \varkappa_j \right\| \right.$$
$$\left. + \sup_{[\overline{T}\lambda_0/\pi] \le k \le \overline{T}} \left\| \frac{(1 - k/\widetilde{T})^{-\delta/2}}{\widetilde{T}^{1/2}} \sum_{j=k+1}^{\widetilde{T}} \gamma_{\theta_0,j} \eta_j \right\| \right\},$$

for any $0 < \delta < 1$. The first factor on the right-hand side of (43) is bounded by

$$K \left| \frac{1}{\widetilde{\widetilde{T}}} \sum_{k=[\overline{T}\lambda_0/\pi]+1}^{\overline{T}} \|\gamma_{\theta_0,k}\| \left(1 - \frac{k}{\widetilde{\widetilde{T}}}\right)^{\delta/2-1} \right| \le K \left( \frac{\overline{T} - [\overline{T}\lambda_0/\pi]}{\widetilde{T}} \right)^{\delta/2},$$

using

$$\|A_{\theta_0,T}^{-1}(k)\| \le K \left(1 - \frac{k}{\widetilde{T}}\right)^{-1},$$

because $\|A_{\theta_0}(\lambda)\| \ge K^{-1}(\pi - \lambda)$ by assumption A5 and because Lemma 1 implies that $\sup_{[\overline{T}\lambda_0/\pi] \le k \le \overline{T}} \|A_{\theta_0,T}(k) - A_{\theta_0}([k\pi/\widetilde{T}])\| = O(T^{-1}\log^2 T)$.

Next, by Lemma 9 the second term inside the braces on the right-hand side of (43) is $O_p(1)$ for $\delta > 0$ small enough, whereas Lemma 8 and (35) imply that the first term is bounded by

$$\begin{aligned}
&\sup_{[\overline{T}\lambda_0/\pi]\le k\le\overline{T}}\left\|\frac{(1-k/\widetilde{T})^{-\delta/2}}{\widetilde{T}}\sum_{j=k+1}^{\widetilde{T}}\gamma_{\theta_0,j}\phi'_{\theta_0,j}\right\|O_p(1)\\
&\quad+O_p\biggl(\sup_{[\overline{T}\lambda_0/\pi]\le k\le\overline{T}}\frac{(1-k/\widetilde{T})^{-\delta/2}}{T^{\delta}}\biggr)\\
&\quad=O_p(|\pi-\lambda_0|^{\delta/2}),
\end{aligned}$$

because of $T^{-1}\le\widetilde{T}^{-1}\le\inf_{[\overline{T}\lambda_0/\pi]\le k\le\overline{T}}(1-k/\widetilde{T})$, $0<\delta<1$, and an obvious extension of Lemma 1 but with $\zeta(\lambda)=\gamma_{\theta_0}(\lambda)\phi'_{\theta_0}(\lambda)$ there. So, (43) is $O_p(|\pi-\lambda_0|^{\delta})$, which implies that (42) holds because $\delta>0$. □

Lemma 11. *Assuming* A1–A6,

$$\sup_{\lambda\in[0,\pi]}\left\|\frac{1}{\widetilde{T}^{1/2}}\sum_{j=[\widetilde{T}\lambda/\pi]+1}^{\widetilde{T}}(\phi_{\theta_T,j}-\phi_{\theta_0,j})\biggl(\frac{I_{X,j}}{h_{\theta_T,j}}-\frac{\sigma^2}{2\pi}\biggr)\right\|=O_p\biggl(\frac{\log T}{T^{1/2}}\biggr). \tag{44}$$

Proof. The expression inside the norm on the left-hand side of (44) is

$$\begin{aligned}
&\frac{1}{\widetilde{T}^{1/2}}\sum_{j=[\widetilde{T}\lambda/\pi]+1}^{\widetilde{T}}\dot{\phi}_{\theta_0,j}\biggl(\frac{I_{X,j}}{h_{\theta_T,j}}-I_{\varepsilon,j}\biggr)(\theta_T-\theta_0)\\
&\qquad+\frac{1}{\widetilde{T}^{1/2}}\sum_{j=[\widetilde{T}\lambda/\pi]+1}^{\widetilde{T}}\dot{\phi}_{\theta_0,j}\biggl(I_{\varepsilon,j}-\frac{\sigma^2}{2\pi}\biggr)(\theta_T-\theta_0) \qquad (45)\\
&\qquad+\frac{1}{\widetilde{T}^{1/2}}\sum_{j=[\widetilde{T}\lambda/\pi]+1}^{\widetilde{T}}(\phi_{\theta_T,j}-\phi_{\theta_0,j}-\dot{\phi}_{\theta_0,j}(\theta_T-\theta_0))\biggl(\frac{I_{X,j}}{h_{\theta_T,j}}-\frac{\sigma^2}{2\pi}\biggr).
\end{aligned}$$

By A6 and then noting that $|a-b|\le(a-b)+2b$ for $a>0$ and $b>0$, the norm of the third term of (45) is bounded by

$$\begin{aligned}
&K\frac{\|\theta_T-\theta_0\|^2}{\widetilde{T}^{1/2}}\sum_{j=1}^{\widetilde{T}}|\log(\lambda_j)|\left|\frac{I_{X,j}}{h_{\theta_T,j}}-\frac{\sigma^2}{2\pi}\right|\\
&\qquad\le K\frac{\|\theta_T-\theta_0\|^2}{\widetilde{T}^{1/2}}\left\{\sum_{j=1}^{\widetilde{T}}|\log(\lambda_j)|\biggl(\frac{I_{X,j}}{h_{\theta_T,j}}-\frac{\sigma^2}{2\pi}\biggr)+\frac{\sigma^2}{\pi}\sum_{j=1}^{\widetilde{T}}|\log\lambda_j|\right\}\\
&\qquad=O_p\biggl(\frac{\log T}{T^{1/2}}\biggr)
\end{aligned}$$

by (35) and then using Lemmas 8 and 7 with $\zeta(\lambda) = |\log\lambda|$, and Lemma 1, respectively. So, uniformly in $\lambda$, the third term of (45) is $o_p(1)$. Likewise, the first term of (45) is $O_p(T^{-1/2})$ uniformly in $\lambda$ using Lemma 8 with $\zeta(\lambda) = \dot{\phi}_{\theta_0}(\lambda)$ and (35). Observe that $\dot{\phi}_{\theta_0}(\lambda)$ satisfies the same conditions as $\zeta(\lambda)$ in Lemma 8 by A6. Finally, the second term of (45) is $O_p(T^{-1/2})$ by Lemma 7 with $\zeta(\lambda) = \dot{\phi}_{\theta_0}(\lambda)$. □

LEMMA 12. *Assuming* A1–A6, *for all* $\epsilon > 0$,

$$\lim_{\lambda_0\to\pi}\limsup_{T\to\infty}\Pr\left\{\sup_{\lambda_0\le\lambda\le\pi}\left|\frac{1}{\widetilde{T}}\sum_{q=[\overline{T}\lambda_0/\pi]+1}^{[\overline{T}\lambda/\pi]}\frac{H_{\theta_T,T}(q)}{\widetilde{T}^{1/2}}\right.\right.$$
$$\left.\left.\times\sum_{j=q+1}^{\widetilde{T}}\gamma_{\theta_T,j}\left(\frac{I_{X,j}}{h_{\theta_T,j}}-\frac{\sigma^2}{2\pi}\right)\right|>\epsilon\right\}=0. \tag{46}$$

PROOF. Notice that (35) implies that it suffices to show (46) in the set $\{\|\theta_T - \theta_0\| < KT^{-1/2}m_T^{-1}\}$, where $m_T + m_T^{-1}T^{-1/2} \to 0$. On the other hand, Lemma 11 and then Lemma 8 imply that, uniformly in $q$,

$$\frac{1}{\widetilde{T}^{1/2}}\sum_{j=q+1}^{\widetilde{T}}\gamma_{\theta_T,j}\varkappa_j = \left(\frac{\sigma^2}{\widetilde{T}}\sum_{j=q+1}^{\widetilde{T}}\gamma_{\theta_0,j}\phi'_{\theta_0,j}\right)\widetilde{T}^{1/2}(\theta_0-\theta_T)+O_p(T^{-\delta}),$$
$$\frac{1}{\widetilde{T}^{1/2}}\sum_{j=q+1}^{\widetilde{T}}\gamma_{\theta_T,j}\eta_j = \frac{1}{\widetilde{T}^{1/2}}\sum_{j=q+1}^{\widetilde{T}}\gamma_{\theta_0,j}\eta_j + O_p(T^{-1/2}), \tag{47}$$

proceeding as in the proof of (44) but with $\varkappa_j + \eta_j$ replaced by $\eta_j$ there. Observe that we can take $\lambda_0 > \pi/2$. Next, uniformly in $q$, A6 implies that

$$\sup_{[\overline{T}\lambda_0/\pi]\le q\le\overline{T}}\|A_{\theta_T,T}(q) - A_{\theta_0,T}(q)\| = (\pi-\lambda_0)O_p(\|\theta_T-\theta_0\|),$$

which will imply that, with probability approaching one, as $T\to\infty$,

$$\|A^{-1}_{\theta_T,T}(q)\| \le \|A^{-1}_{\theta_0,T}(q)\|(1+KT^{-1/2}m_T^{-1}) \le K\left(1-\frac{q}{\widetilde{T}}\right)^{-1},$$

because $\|A_{\theta_0}(\lambda)\| \ge K^{-1}(\pi-\lambda)$ and Lemma 1 implies that

$$\sup_{[\overline{T}\lambda_0/\pi]\le q\le\overline{T}}\|A_{\theta_0,T}(q) - A_{\theta_0}([q\pi/\widetilde{T}])\| = O(T^{-1}\log^2 T).$$

So, we have that, for $0 < \delta < 1/2$,

$$\sup_{\lambda_0\le\lambda\le\pi}\left\|\frac{1}{\widetilde{T}}\sum_{q=[\overline{T}\lambda_0/\pi]+1}^{[\overline{T}\lambda/\pi]}\frac{H_{\theta_T,T}(q)}{\widetilde{T}^{1/2}}\sum_{j=q+1}^{\widetilde{T}}\gamma_{\theta_T,j}\left(\frac{I_{X,j}}{h_{\theta_T,j}}-\frac{\sigma^2}{2\pi}\right)\right\|$$

$$
\begin{aligned}
&\le K \sup_{\lambda_0\le\lambda\le\pi}\left|\frac{1}{\widetilde{T}}\sum_{q=[\overline{T}\lambda_0/\pi]+1}^{[\overline{T}\lambda/\pi]}\|\gamma_{\theta_0,q}\|\left(1-\frac{q}{\widetilde{T}}\right)^{-1+\delta/2}\right| \\
&\qquad\times\left\{\sup_{[\overline{T}\lambda_0/\pi]\le q\le\overline{T}}\left\|\left(1-\frac{q}{\widetilde{T}}\right)^{-\delta/2}\frac{1}{\widetilde{T}^{1/2}}\sum_{j=q+1}^{\widetilde{T}}\gamma_{\theta_0,j}\eta_j\right\|\right. \\
&\qquad\qquad\qquad\left.+O_p(|\pi-\lambda_0|^{\delta/2})\right\},
\end{aligned}
\tag{48}
$$

by (47) and because $T^{-1}\le\widetilde{T}^{-1}\le\inf_{[\overline{T}\lambda_0/\pi]\le q\le\overline{T}}(1-q/\widetilde{T})$. But Lemma 9 implies that

$$
\sup_{[\overline{T}\lambda_0/\pi]\le q\le\overline{T}}\left\|(1-q/\widetilde{T})^{-\delta/2}\widetilde{T}^{-1/2}\sum_{j=q+1}^{\widetilde{T}}\gamma_{\theta_0,j}\eta_j\right\|=O_p(1),
$$

and A3 implies that

$$
\sup_{\lambda_0\le\lambda\le\pi}\frac{1}{\widetilde{T}}\sum_{q=[\overline{T}\lambda_0/\pi]+1}^{[\overline{T}\lambda/\pi]}\|\gamma_{\theta_0,q}\|\left(1-\frac{q}{\widetilde{T}}\right)^{-1+\delta/2}\le K\left(\frac{\overline{T}-[\overline{T}\lambda_0/\pi]}{\widetilde{T}}\right)^{\delta/2},
$$

and, hence, the left-hand side of (48) is $O_p(|\pi-\lambda_0|^{\delta/2})$. From here we conclude that (46) holds because $\delta>0$. □

**7. Proofs.** This section provides the proofs of the main results which are based on the series of lemmas given in the previous section.

Proof of Theorem 1. Part (a) follows by Lemma 4 with $\zeta(\lambda)=1$ there. The proof of part (b) follows immediately from part (a) and Lemma 7 with $\zeta(\lambda)=1$ there. □

Proof of Theorem 2. Part (a). By Lemma 8 with $\zeta(\lambda)=1$ there and the definitions of $G_{\theta,T}$ and $G_T^0$, we have that

$$
\begin{aligned}
&\widetilde{T}^{1/2}(G_{\theta_T,T}(\lambda)-G_T^0(\lambda)) \\
&\qquad=-\left(\frac{\sigma^2}{\widetilde{T}}\sum_{j=1}^{[\widetilde{T}\lambda/\pi]}\phi'_{\theta_0,j}\right)\widetilde{T}^{1/2}(\theta_T-\theta_0)+o_p(1) \\
&\qquad=-\left(\frac{\sigma^2}{\widetilde{T}}\sum_{j=1}^{[\widetilde{T}\lambda/\pi]}\phi'_{\theta_0,j}\right)S_T^{-1}\frac{2\pi}{G_{\theta_0,T}(\pi)\widetilde{T}^{1/2}}\sum_{k=1}^{\widetilde{T}}\phi_{\theta_0,k}\frac{I_{X,k}}{h_{\theta_0,k}} \\
&\qquad\quad+o_p(1),
\end{aligned}
\tag{49}
$$

by (8) and (9), and where the $o_p(1)$ term is uniform in $\lambda \in [0,\pi]$. Likewise,

$$\widetilde{T}^{1/2}(G_{\theta_T,T}(\pi) - G_T^0(\pi)) = o_p(1) \tag{50}$$

because of (36) and (37) and, by Lemma 1 with $\zeta(\lambda) = \phi_{\theta_0}(\lambda)$ and (9), we have $\|\widetilde{T}^{-1}\sum_{j=1}^{\widetilde{T}} \phi_{\theta_0,j}\| = O(T^{-1}\log T)$. So, (50) holds. Also, it is worth noticing that Lemma 1 with $\zeta(\lambda) = \phi_{\theta_0}(\lambda)\phi'_{\theta_0}(\lambda)$ implies that $\|S_T - \Sigma_{\theta_0}\| = O(T^{-1}\log^2 T)$.

On the other hand, noting that (50) and A1 imply that

$$G_T^0(\pi) = \sigma^2 + O_p(T^{-1/2}), \tag{51}$$

and that

$$|G_{\theta_0,T}(\pi) - G_T^0(\pi)| = o_p(\widetilde{T}^{-1/2})$$

by Lemma 4, then by (49), (50) and (36), uniformly in $\lambda$, we obtain that

$$\begin{aligned}\alpha_{\theta_T,T}(\lambda) &= \alpha_T^0(\lambda) + \frac{\widetilde{T}^{1/2}(G_{\theta_T,T}(\lambda) - G_T^0(\lambda))}{G_T^0(\pi)} \\ &\quad + G_{\theta_T,T}(\lambda)\widetilde{T}^{1/2}\left(\frac{1}{G_{\theta_T,T}(\pi)} - \frac{1}{G_T^0(\pi)}\right) \\ &= \alpha_T^0(\lambda) - \frac{1}{\widetilde{T}}\sum_{j=1}^{[\widetilde{T}\lambda/\pi]}\left[\phi'_{\theta_0,j}S_T^{-1}\frac{2\pi}{G_T^0(\pi)\widetilde{T}^{1/2}}\sum_{k=1}^{\widetilde{T}}\phi_{\theta_0,k}I_{\varepsilon,k}\right] + o_p(1),\end{aligned} \tag{52}$$

which concludes the proof of part (a).

Next part (b). Taking into account part (a), part (b) follows because Lemma 7 guarantees the fidi's convergence of $\alpha_T^0$ and its tightness. Tightness of the second term on the right-hand side of (52) follows by (37) and Lemma 1 and then because $\int_0^\lambda \phi_{\theta_0}(u)\,du$ is Hölder's continuous of order greater than 1/2 by A3. This concludes the proof of the theorem. □

PROOF OF THEOREM 3. Using (51) and recalling that $H_{\theta,T}(j) = \gamma'_{\theta,j}A_{\theta,T}^{-1}(j)$, we obtain that

$$\begin{aligned}\beta_T^0(\lambda) &= \frac{1}{\widetilde{T}^{1/2}}\sum_{j=1}^{[\widetilde{T}\lambda/\pi]}\left(\left(\frac{2\pi}{\sigma^2}I_{\varepsilon,j} - 1\right)\right. \\ &\qquad \left. - H_{\theta_0,T}(j)\frac{1}{\widetilde{T}}\sum_{k=j+1}^{\widetilde{T}}\gamma_{\theta_0,k}\left(\frac{2\pi}{\sigma^2}I_{\varepsilon,k} - 1\right)\right) + o_p(1),\end{aligned} \tag{53}$$

where the $o_p(1)$ term is uniform in $\lambda \in [0,\pi]$.

Suppose, to be shown later, that the convergence in $[0,\lambda_0]$ holds for any $0<\lambda_0<\pi$. Then, because $B_\pi$ and the limit of the process $\widetilde{T}^{-1/2}\sum_{j=1}^{[\widetilde{T}\lambda/\pi]}(I_{\varepsilon,j}-\sigma^2/2\pi)$ are continuous in $[0,\pi]$, Billingsley's [3] Theorem 4.2 implies that it suffices to show that, for all $\epsilon>0$,

$$\lim_{\lambda_0\to\pi}\limsup_{T\to\infty}\Pr\left\{\sup_{\lambda_0\le\lambda\le\pi}\left|\frac{1}{\widetilde{T}}\sum_{j=[\widetilde{T}\lambda_0/\pi]+1}^{[\widetilde{T}\lambda/\pi]}\frac{H_{\theta_0,T}(j)}{\widetilde{T}^{1/2}}\right.\right.$$

$$\left.\left.\times\sum_{k=j+1}^{\widetilde{T}}\gamma_{\theta_0,k}\left(\frac{2\pi}{\sigma^2}I_{\varepsilon,k}-1\right)\right|>\epsilon\right\}=0,$$

which follows by Lemma 10; compare the second term on the right-hand side of (43).

So, to complete the proof, we need to show that, for any $0<\lambda_0<\pi$,

$$\begin{aligned}&\frac{1}{\widetilde{T}^{1/2}}\sum_{j=1}^{[\widetilde{T}\lambda/\pi]}\left(\left(\frac{2\pi}{\sigma^2}I_{\varepsilon,j}-1\right)\right.\\ &\qquad\left.-H_{\theta_0,T}(j)\frac{1}{\widetilde{T}}\sum_{k=j+1}^{\widetilde{T}}\gamma_{\theta_0,k}\left(\frac{2\pi}{\sigma^2}I_{\varepsilon,k}-1\right)\right)\Rightarrow\frac{1}{\pi^{1/2}}B_\pi(\lambda),\end{aligned}\tag{54}$$

in $[0,\lambda_0]$. Fidi's convergence follows by Lemma 7, part (b) after we note that the second term on the right-hand side of (53) is

$$\frac{1}{\widetilde{T}^{1/2}}\sum_{k=1}^{\widetilde{T}}\left(\frac{1}{\widetilde{T}}\sum_{j=1}^{k\wedge[\widetilde{T}\lambda/\pi]}H_{\theta_0,T}(j)\right)\gamma_{\theta_0,k}\left(\frac{2\pi}{\sigma^2}I_{\varepsilon,k}-1\right)$$

and $(\widetilde{T}^{-1}\sum_{j=1}^{k\wedge[\widetilde{T}\lambda/\pi]}H_{\theta_0,T}(j))\gamma_{\theta_0,k}$ satisfies the same conditions of Lemma 7 for $\zeta(\lambda)$, for example, those of $h_n(\lambda)$ in [14], Theorem 4.2. Then, it suffices to prove tightness. Since $\alpha_T^0$ is tight, we only need to show the tightness condition of

$$\Lambda_T(\lambda)=\frac{1}{\widetilde{T}}\sum_{j=1}^{[\widetilde{T}\lambda/\pi]}H_{\theta_0,T}(j)\left(\frac{1}{\widetilde{T}^{1/2}}\sum_{k=j+1}^{\widetilde{T}}\gamma_{\theta_0,k}\left(I_{\varepsilon,k}-\frac{\sigma^2}{2\pi}\right)\right).\tag{55}$$

By Billingsley's [3] Theorem 15.6, it suffices to show that

$$\mathbb{E}(|\Lambda_T(\vartheta)-\Lambda_T(\mu)||\Lambda_T(\lambda)-\Lambda_T(\vartheta)|)\le K|\lambda-\mu|^{2\delta}$$

for all $0\le\mu<\vartheta<\lambda\le\pi$ and some $\delta>1/2$. Observe that we can take $\widetilde{T}^{-1}<|\lambda-\mu|$, since otherwise the last inequality is trivial. Because $(\lambda-\vartheta)(\vartheta-\mu)<$

$(\lambda-\mu)^2$, by the Cauchy–Schwarz inequality, it suffices to show the last displayed inequality holds for $\mathbb{E}|\Lambda_T(\lambda)-\Lambda_T(\mu)|^2$, which is

$$\begin{aligned}
&\frac{1}{\widetilde{T}^3}\sum_{j,k=[\overline{T}\mu/\pi]+1}^{[\overline{T}\lambda/\pi]} H_{\theta_0,T}(j)\sum_{\ell_1=j+1}^{\widetilde{T}}\sum_{\ell_2=k+1}^{\widetilde{T}}\gamma_{\theta_0,\ell_1}\gamma'_{\theta_0,\ell_2}\\
&\qquad\qquad\times\mathbb{E}\biggl[\biggl(I_{\varepsilon,\ell_1}-\frac{\sigma^2}{2\pi}\biggr)\biggl(I_{\varepsilon,\ell_2}-\frac{\sigma^2}{2\pi}\biggr)\biggr]H'_{\theta_0,T}(k)\\
&\quad\le\frac{K}{\widetilde{T}^2}\sum_{j,k=[\overline{T}\mu/\pi]+1}^{[\overline{T}\lambda/\pi]}\|H_{\theta_0,T}(j)\|\|H_{\theta_0,T}(k)\|\\
&\quad\le K(|\widetilde{H}(\lambda)-\widetilde{H}(\mu)|^2+\widetilde{T}^{-2}\log^2\widetilde{T}),
\end{aligned}$$

where

$$\widetilde{H}(\lambda):=\pi^{-1}\int_0^\lambda H_{\theta_0}(x)\,dx\quad\text{and}\quad\|\widetilde{H}_T(\lambda)-\widetilde{H}(\lambda)\|\le K\widetilde{T}^{-1}\log T$$

and

$$\widetilde{H}_T(\lambda):=\widetilde{T}^{-1}\sum_{j=1}^{[\overline{T}\lambda/\pi]}\|H_{\theta_0,T}(j)\|$$

by Lemma 1. From here we conclude by Billingsley's [3] Theorem 15.6, because $\widetilde{H}(\lambda)$ is a monotonic, continuous and nondecreasing function such that $|\widetilde{H}(\lambda)-\widetilde{H}(\mu)|\le K|\lambda-\mu|^\delta$, $\delta>1/2$ and $\widetilde{T}^{-1}\le|\lambda-\mu|$. □

Proof of Theorem 4. By definition of $\beta_{\theta,T}$ and $\beta_T^0$, it suffices to show that

$$\left|\frac{1}{\widetilde{T}^{1/2}}\sum_{k=1}^{[\overline{T}\lambda/\pi]}\biggl(\frac{I_{X,k}}{h_{\theta_T,k}}-I_{\varepsilon,k}\biggr)-H_{\theta_0,T}(k)\frac{1}{\widetilde{T}}\sum_{j=k+1}^{\widetilde{T}}\gamma_{\theta_0,j}\biggl(\frac{I_{X,j}}{h_{\theta_T,j}}-I_{\varepsilon,j}\biggr)\right| \tag{56}$$

and

$$\begin{aligned}
&\frac{1}{G_{\theta_T,T}(\pi)}\Biggl(\frac{1}{\widetilde{T}}\sum_{k=1}^{[\overline{T}\lambda/\pi]}H_{\theta_0,T}(k)\frac{1}{\widetilde{T}^{1/2}}\sum_{j=k+1}^{\widetilde{T}}\gamma_{\theta_0,j}\biggl(\frac{I_{X,j}}{h_{\theta_T,j}}-\frac{G_{\theta_T,T}(\pi)}{2\pi}\biggr)\Biggr)\\
&\qquad-\frac{1}{G_{\theta_T,T}(\pi)}\\
&\qquad\quad\times\Biggl(\frac{1}{\widetilde{T}}\sum_{k=1}^{[\overline{T}\lambda/\pi]}H_{\theta_T,T}(k)\frac{1}{\widetilde{T}^{1/2}}\sum_{j=k+1}^{\widetilde{T}}\gamma_{\theta_T,j}\biggl(\frac{I_{X,j}}{h_{\theta_T,j}}-\frac{G_{\theta_T,T}(\pi)}{2\pi}\biggr)\Biggr)
\end{aligned} \tag{57}$$

converge to zero uniformly in $\lambda \in [0,\pi]$. Expression (56) is $o_p(1)$, uniformly in $\lambda \in [0,\pi]$, because the contribution due to the term in brackets in the last line of (52), that is, $-\phi'_{\theta_0,j} 2\pi (G_T^0(\pi))^{-1} S_T^{-1} \widetilde{T}^{-1/2} \sum_{k=1}^{\widetilde{T}} \phi_{\theta_0,k} I_{\varepsilon,k}$, is easily seen to be zero. Next, because

$$\begin{aligned} &\frac{1}{\widetilde{T}} \sum_{k=1}^{[\widetilde{T}\lambda/\pi]} \|\gamma_{\theta_0,k}\| \|A_{\theta_0,T}^{-1}(k)\| \frac{1}{\widetilde{T}} \sum_{j=k+1}^{\widetilde{T}} \|\gamma_{\theta_0,j}\| \\ &\qquad \le K \frac{1}{\widetilde{T}} \sum_{k=1}^{[\widetilde{T}\lambda/\pi]} \|\gamma_{\theta_0,k}\| \left\| A_{\theta_0,T}^{-1}(k) \left(1 - \frac{k}{\widetilde{T}}\right) \right\| \\ &\qquad \le K \frac{1}{\widetilde{T}} \sum_{k=1}^{[\widetilde{T}\lambda/\pi]} \|\gamma_{\theta_0,k}\| \le K \end{aligned}$$

by integrability of $\gamma_{\theta_0}$ and $\|A_{\theta_0,T}(k)(1-k/\widetilde{T})^{-1}\| > 0$ by A3 and A5, it implies that the contribution to (56) due to the term $o_p(1)$ on the right-hand side of (52) is negligible.

Next we examine (57). Because of (50) and (51), it suffices to show that

$$\begin{aligned} &\frac{1}{\widetilde{T}} \sum_{k=1}^{[\widetilde{T}\lambda/\pi]} \left\{ \frac{H_{\theta_0,T}(k)}{\widetilde{T}^{1/2}} \sum_{j=k+1}^{\widetilde{T}} \gamma_{\theta_0,j} \left( \frac{I_{X,j}}{h_{\theta_T,j}} - \frac{\sigma^2}{2\pi} \right) \right. \\ &\qquad \left. - \frac{H_{\theta_T,T}(k)}{\widetilde{T}^{1/2}} \sum_{j=k+1}^{\widetilde{T}} \gamma_{\theta_T,j} \left( \frac{I_{X,j}}{h_{\theta_T,j}} - \frac{\sigma^2}{2\pi} \right) \right\} \end{aligned} \tag{58}$$

converges to zero uniformly in $\lambda \in [0,\pi]$, after observing that

$$\sup_{\lambda\in[0,\pi]} \left| \sum_{k=1}^{[\widetilde{T}\lambda/\pi]} H_{\theta_T,T}(k) \sum_{j=k+1}^{\widetilde{T}} \gamma_{\theta_T,j} - \sum_{k=1}^{[\widetilde{T}\lambda/\pi]} H_{\theta_0,T}(k) \sum_{j=k+1}^{\widetilde{T}} \gamma_{\theta_0,j} \right| = 0.$$

First, we observe that Lemmas 10 and 12 imply that it suffices to show the uniform convergence in $\lambda \in [0,\lambda_0]$ for any $\lambda_0 < \pi$. But (58) is equal to

$$\frac{1}{\widetilde{T}} \sum_{k=1}^{[\widetilde{T}\lambda/\pi]} H_{\theta_T,T}(k) \frac{1}{\widetilde{T}^{1/2}} \sum_{j=k+1}^{\widetilde{T}} (\gamma_{\theta_0,j} - \gamma_{\theta_T,j}) \left( \frac{I_{X,j}}{h_{\theta_T,j}} - \frac{\sigma^2}{2\pi} \right) \tag{59}$$

$$+ \frac{1}{\widetilde{T}} \sum_{k=1}^{[\widetilde{T}\lambda/\pi]} (H_{\theta_0,T}(k) - H_{\theta_T,T}(k)) \frac{1}{\widetilde{T}^{1/2}} \sum_{j=k+1}^{\widetilde{T}} \gamma_{\theta_0,j} \left( \frac{I_{X,j}}{h_{\theta_T,j}} - \frac{\sigma^2}{2\pi} \right). \tag{60}$$

So, the theorem follows if (59) and (60) are $o_p(1)$ uniformly in $\lambda \in [0,\lambda_0]$.

To that end, we first show that

$$\sup_{\lambda\in[0,\pi]} \frac{1}{\widetilde{T}} \sum_{j=1}^{[\widetilde{T}\lambda/\pi]} \|\phi_{\theta_0,j} - \phi_{\theta_T,j}\| = o_p(1), \tag{61}$$

$$\sup_{\lambda\in[0,\lambda_0]} \|A^{-1}_{\theta_0,T}(\lambda) - A^{-1}_{\theta_0}(\lambda)\| = o(1), \tag{62}$$

$$\sup_{\lambda\in[0,\lambda_0]} \|A^{-1}_{\theta_T,T}(\lambda) - A^{-1}_{\theta_0,T}(\lambda)\| = o_p(1). \tag{63}$$

(61) follows proceeding as with the proof of (44) in Lemma 11, but without the factor $h^{-1}_{\theta_T,j} I_{X,j} - \sigma^2/(2\pi)$; (62) follows because assumption A5 implies that $A_{\theta_0}(\lambda_0) > 0$ and because, by assumption A3, $\|\phi_{\theta_0}(\lambda)\phi'_{\theta_0}(\lambda)\|$ satisfies the same conditions on $\zeta(\lambda)$ in Lemma 1, so that

$$\sup_{\lambda\in[0,\lambda_0]} \|A_{\theta_0}(\lambda) - A_{\theta_0,T}(\lambda)\| = O(T^{-1}\log^2 T);$$

and (63) follows proceeding as with the proof of (61) and (62).

Now we show that (59) is $o_p(1)$ uniformly in $\lambda \in [0, \lambda_0]$, which follows by Lemma 11 and (61)–(63), noting that $(\gamma'_{\theta_0,j} - \gamma'_{\theta_T,j}) = (0, \phi'_{\theta_0,j} - \phi'_{\theta_T,j})$; so does (60) by (61) and (63) and noting that

$$\sup_{\lambda\in[0,\pi]} \left| \frac{1}{\widetilde{T}^{1/2}} \sum_{j=[\widetilde{T}\lambda/\pi]+1}^{\widetilde{T}} \gamma_{\theta_0,j} \left( \frac{I_{X,j}}{h_{\theta_T,j}} - \frac{\sigma^2}{2\pi} \right) \right| = O_p(1)$$

by Lemmas 7 and 8 with $\zeta(\lambda) = \gamma_{\theta_0}(\lambda)$ there and observing (35) and that by Lemma 1, $\widetilde{T}^{-1} \sum_{j=[\widetilde{T}\lambda/\pi]+1}^{\widetilde{T}} \gamma_{\theta_0,j}\phi'_{\theta_0,j} \to \int_\lambda^\pi \gamma_{\theta_0}(x)\phi'_{\theta_0}(x)\,dx$. □

Proof of Theorem 5. Under $H_{1T}$, we have that, by definition,

$$G_{\theta_0,T}(\lambda) = \frac{2\pi}{\widetilde{T}} \sum_{j=1}^{[\widetilde{T}\lambda/\pi]} \frac{I_{X,j}}{h_j} + \frac{\sigma^2\tau}{\widetilde{T}^{3/2}} \sum_{j=1}^{[\widetilde{T}\lambda/\pi]} l_j$$
$$+ \frac{2\pi\tau}{\widetilde{T}^{3/2}} \sum_{j=1}^{[\widetilde{T}\lambda/\pi]} l_j \left( \frac{I_{X,j}}{h_j} - \frac{\sigma^2}{2\pi} \right) + \frac{1}{\widetilde{T}^2} \sum_{j=1}^{[\widetilde{T}\lambda/\pi]} s_{T,j} \frac{I_{X,j}}{h_j}.$$

By Lemmas 1, 4 and 7 with $\zeta(\lambda) = \tau l(\lambda)$, and because $|s_T|$ is integrable, we have

$$G_{\theta_0,T}(\lambda) = \frac{2\pi}{\widetilde{T}} \sum_{j=1}^{[\widetilde{T}\lambda/\pi]} I_{\varepsilon,j} + \frac{\sigma^2\tau}{\widetilde{T}^{1/2}\pi} \int_0^\lambda l(u)\,du + o_p(T^{-1/2}).$$

So, using (51) because $\int_0^\pi l(u)\,du = 0$, we have that, uniformly in $\lambda \in [0, \pi]$,

$$\widetilde{T}^{1/2}\left(\frac{G_{\theta_0,T}(\lambda)}{G_{\theta_0,T}(\pi)} - \frac{\lambda}{\pi}\right) = \widetilde{T}^{1/2}\left(\frac{2\pi}{G_T^0(\pi)\widetilde{T}} \sum_{j=1}^{[\widetilde{T}\lambda/\pi]} I_{\varepsilon,j} - \frac{\lambda}{\pi} + \frac{\tau}{\widetilde{T}^{1/2}\pi}\int_0^\lambda l(u)\,du\right) + o_p(1)$$
$$= \alpha_T^0(\lambda) + \frac{\tau}{\pi}\int_0^\lambda l(u)\,du + o_p(1).$$

From here the conclusion is straightforward. □

**Acknowledgments.** We thank an Associate Editor, two referees and Hira Koul for their constructive comments on previous versions of this article which have led to substantial improvement of the paper. Of course, all remaining errors are our sole responsibility.

M. A. Delgado
C. Velasco
Departamento de Economía
Universidad Carlos III de Madrid
C./Madrid 126-128
Getafe, 28903 Madrid
Spain
e-mail: miguelangel.delgado@uc3m.es
carlos.velasco@uc3m.es

J. Hidalgo
Department of Economics
London School of Economics
Houghton Street
London W2A 2AE
United Kingdom
e-mail: f.j.hidalgo@lse.ac.uk